\documentclass[proc]{edpsmath}
\usepackage{xcolor}
\usepackage{etoolbox}
\usepackage{graphicx}

\newtheorem{thm}{Theorem}

\newtheorem{rmk}{Remark}

\newcommand{\indic}[1]{\mathbf{1}_{#1}}
\newcommand{\Rb}{\ensuremath{\mathbb R}}
\newcommand{\Nb}{\ensuremath{\mathbb N}}

\begin{document}

%%-----------------------------
%%      the top matter
%%-----------------------------
\title{Modeling Compartmentalization within Intracellular Signaling Pathway}\thanks{We acknowledge the funding from the IMAGO project, Digit-Bio INRAE Metaprogram}\thanks{This work is supported by a public grant overseen by the French National research Agency (ANR) as part of the « Investissements d’Avenir » program, through the "ADI 2020" project funded by the IDEX Paris-Saclay, ANR-11-IDEX-0003-02}\thanks{This work is supported by the Italian Ministry of Education, Universities and Research through the MIUR grant “Dipartimenti di Eccellenza 2018-2022” (Project no. E11G18000350001)}\thanks{This work is partially supported by the spanish  MINECO-Feder (grant RTI2018-098850-B-I00) and Junta de Andaluc\'ia (grants PY18-RT-2422 and  A-FQM-311-UGR18)}
\thanks{...}% At most 5 thanks
\author{Claire Alamichel}\address{Université Paris-Saclay, CNRS UMR 8071, Univ Evry, Laboratoire de Mathématiques et Modélisation d'Evry, 91037, Evry-Courcouronnes, France \& Department of Mathematical Sciences, Politecnico di Torino, Corso Duca degli Abruzzi 24, Torino, 10129, Italy.}
\author{Juan Calvo}\address{Universidad de Granada, Avda. Hospicio s/n, 18071, Granada, Spain \&  Research Unit ``Modeling Nature'' (MNat).}
\author{Erwan Hingant}\address{Université de Picardie Jules Verne, LAMFA, CNRS UMR 7352, 33 rue Saint-Leu, 80039, Amiens, France.}
\author{Saoussen Latrach}\address{Laboratoire Analyse, Géométrie et Application, Université Sorbonne Paris Nord.}
\author{Nathan Quiblier}\address{AIstroSight, Inria, Université Claude Bernard Lyon 1, Hospices Civils de Lyon \& Liris, UMR 5205 CNRS, INSA Lyon, Université Claude Bernard Lyon 1.}
\author{Romain Yvinec}\address{Université Paris-Saclay, Inria, Inria Saclay-Île-de-France, 91120, Palaiseau, France\\
PRC, INRAE, CNRS, Université de Tours, 37380 Nouzilly, France.}
%
%\dedicated{\it Dedicated to Maurice Dupont} %if necessary
%
\begin{abstract} 

We present a new modeling approach for G protein coupled receptors signaling systems, that take into account the compartmentalization of receptors and their effectors, both at plasma membrane and in dynamic intra-cellular vesicles called endosomes. The first building block of the model is about compartment dynamics. It takes into account creation of de-novo endosomes, \textit{i.e.} endocytosis, recycling of endosomes back to plasma membrane, degradation through transfer into lysosomes as well as endosomes fusion through coagulation dynamics. The second building block is biochemical reactions into each compartments and the transfer of molecules between the dynamical compartments. In this work, we prove sufficient conditions to obtain exponentially ergodicity for the size distribution of intracellular compartments. We futher design a finite volume scheme to simulate our model and show two application cases for receptor trafficking and spatially biased second effector signaling.

\end{abstract}
\begin{resume} 

Nous présentons une nouvelle approche de modélisation des systèmes de signalisation des récepteurs couplés aux protéines G, qui prend en compte la compartimentation des récepteurs et de leurs effecteurs, à la fois au niveau de la membrane plasmique et dans des vésicules intracellulaires très dynamiques, appelées endosomes. Le premier bloc du modèle concerne la dynamique des compartiments. Il prend en compte la création d'endosomes \textit{de-novo}, l'endocytose, le recyclage des endosomes vers la membrane plasmique, la dégradation par transfert dans les lysosomes ainsi que la fusion des endosomes par une dynamique de coagulation. Le deuxième bloc du modèle est constitué des réactions biochimiques dans chaque compartiment et du transfert de molécules entre les compartiments dynamiques. Dans ce travail, nous prouvons des conditions suffisantes pour obtenir une ergodicité exponentielle pour la distribution en taille des compartiments intracellulaires. Nous concevons de plus un schéma aux volumes finis pour simuler notre modèle et montrons deux cas d'application pour le trafic des récepteurs et la signalisation spatialement biaisée d'effecteurs.

\end{resume}
\maketitle
%%-----------------------------
%%      your text
%%-----------------------------

%\showthe\columnwidth
\section{Introduction}

G Protein Coupled Receptors (GPCR) are a large class of membrane receptors and form an important class of targeted pharmaceutical agents in many different contexts. Recently, it has been shown that those receptors are pleiotropic, i.e. able to selectively activate different signalling pathways. A complex spatio-temporal encoding of the effector molecules downstream the receptors can be modulated by the stimuli (hormone, biochemical ligand...) and shed light into important cell regulation mechanisms that need to be taken into account to design efficient therapeutic strategies \cite{kenakinBiasedReceptorSignaling2019}.

Key evidence have further shown that the intra-cellular traffic of internalized receptors (endocytosis) has a major impact in cell response to a given stimuli \cite{jean-alphonseRegulationGPCRSignal2011,vilardagaEndosomalGenerationCAMP2014}. In particular, for a number of GPCRs, including the Beta-2 Adrenergic Receptors ($\beta$2AR) \cite{kimSpatiotemporalCharacterizationGPCR2021}, the parathyroid hormone receptor
(PTHR) \cite{jean-alphonseSpatiallyRestrictedProteincoupled2014}, the luteinizing hormone/choriogonadotropin receptor (LHCGR) \cite{lygaPersistentCAMPSignaling2016}, or the follicle-stimulating hormone receptor (FSHR) \cite{sayersIntracellularFollicleStimulatingHormone2018}, the production of the second messenger cyclic adenosine monophosphate (cAMP) occurs first at the plasma membrane and then from a highly dynamic pool of intracellular vesicles (called endosomes) following internalization of the receptor by endocytosis. 
This spatio-temporal dimension of signaling has a significant impact on physiological functions, such as the control of serum calcium by PTHR signaling \cite{whiteSpatialBiasCAMP2021}, or the resumption of meiosis by LHCGR signaling \cite{lygaPersistentCAMPSignaling2016}. 

Thus, to faithfully represent the complexity of signalling pathways, we need to take into account the dynamic of the transient pool of specialised intracellular endosomes \cite{sorkinEndocytosisSignallingIntertwining2009,birtwistleEndocytosisSignallingMeeting2009,villasenorSignalProcessingEndosomal2016} following receptor stimulation, and its role on the reaction networks involved in the signalling pathways. The current biological hypothesis is that the endosomal compartments provide a dynamic and an heterogeneous compartmentalised structure that allows specialised molecules to be separated from the bulk cytoplasm (physically separated through a lipid bilayer) in order to have a proper function of the cell response.
    
Classical ways to model the dynamic behaviour of signalling pathways use Chemical Reaction Networks \cite{ingallsMathematicalModelingSystems2013}, either in a deterministic formalism using ordinary differential equations \cite{feinbergFoundationsChemicalReaction2019}, or in a stochastic formalism (typically when few molecules are present) using continuous-time Markov chains \cite{andersonStochasticAnalysisBiochemical2015}. Both approaches typically assume the law of mass action and an idealised homogeneous environment. When spatial dynamics is important to take into account, one may use reaction-diffusion models to represent for instance spatial gradients \cite{kholodenkoCellsignallingDynamicsTime2006}, or compartmental models,  to physically represent segregation between static compartments \cite{weddellIntegrativeMetamodelingIdentifies2017}. Relatively few works have addressed the issue of representing explicitly the segregation of molecules in a dynamic environment. The peculiarity of the compartmentalised signalling pathways is indeed that the (relevant) endosomal compartments are created upon receptor activation, and their number, their size and their molecular content evolved dynamically, somehow within a similar time scale than that of the signalling pathways activation. To the best of our knowledge, a first attempt of defining such models, within the context of signalling pathways, dates back to \cite{foretGeneralTheoreticalFramework2012} and uses a deterministic population dynamics formalism to follow a population of compartments, structured by their size and molecular content, and which undergo coagulation-fragmentation representing endocytosis, fusion, fission, recycling and  degradation. A transport-like operator further represents the chemical reactions within each compartment. Recently, a stochastic counterpart has been proposed by \cite{dusoStochasticReactionNetworks2020,pietzschCompartorToolboxAutomatic2021}.

In this work, we study the long-time behavior and numerical schemes for minimal deterministic models that can represent compartmentalised signalling pathways. We take inspiration from \cite{foretGeneralTheoreticalFramework2012}; our models are able to represent qualitatively main experimental observations from \cite{jean-alphonseSpatiallyRestrictedProteincoupled2014,whiteSpatialBiasCAMP2021}.

In section \ref{sec:models}, we describe two minimal models, structured with respectively one or two variables, that can represent the size-distribution of the endosomal compartment population and their molecular content. In section \ref{sec:longtime}, we study the long time behavior of the model structured with a single variable. In section \ref{sec:numeric}, we present a numerical scheme for the more general model, structured with two variables. In section \ref{sec:application}, we present some numerical simulations that provide qualitative comparisons with experimental observations. 

\section{Modelling compartmentalised signalling pathways}\label{sec:models}

The first objective is to be able to define a well-posed deterministic model to simulate the endosomal compartment dynamics from their size structure perspective only. We therefore adopt a deterministic population dynamic approach, where individuals are structured by a single positive variable (their size). From biological observations, the main processes that shape the size distribution of the endosomal population include:
\begin{itemize}
	\item Endocytosis: creation of a  de-novo compartment from the cell membrane 
	\item Removal of compartment: either recycling back to the cell membrane, or degradation through lysosomal pathways 
	\item Fusion: binary coagulation of compartments 
\end{itemize}
Let $f=f(t,r)$ be the population density of endosomal compartments at time $t$ and size (volume) $r$. The evolution equation for $f$, that takes into account the four mechanisms above is given by:

%\begin{multline}
\begin{equation}
\label{eq:coag_frag}
\tag{M1}
\frac{\partial f}{\partial t}
\!=\!\underbrace{\frac{1}{2}\!\int_0^r \kappa(r-r',r')f(t,r-r')f(t,r')\mathrm{d}r'\!-\!\int_0^\infty \kappa(r,r')f(t,r)f(t,r')\mathrm{d}r'}_{\fcolorbox{red}{white}{coagulation}}
+\underbrace{\alpha(r)}_{\fcolorbox{blue}{white}{endocytosis}}-\underbrace{\gamma(r)f(t,r)}_{\fcolorbox{blue}{white}{removal}}.
%\\
%+\underbrace{\alpha(r)}_{\fcolorbox{blue}{white}{endocytosis}}-\underbrace{\gamma(r)f(t,r)}_{\fcolorbox{blue}{white}{removal}}
\end{equation}
%\end{multline}
A similar equation may be found in \cite{foretGeneralTheoreticalFramework2012,alexandrovDynamicsIntracellularClusters2022}. The coagulation operator with kernel $\kappa$ preserves mass. Endocytosis is a zero-order process at rate $\alpha$ (source term), compartment removal is a first-order process and occurs at a rate $\gamma$.

In the sequel, model \eqref{eq:coag_frag} will be referred as our 1D model.

The second modelling step is to include molecular content into Eq.~\eqref{eq:coag_frag}. One may first think as an additional scalar structure variable\footnote{In future work, this additional structure variable could be an arbitrary finite dimensional variable to represent other molecular actors of signalling pathways that are either physically located at the plasma membrane, in the endosomal compartments or at the vicinity of those.}, for instance the quantity of (active) receptor within each compartment is of primary interest to represent receptor trafficking with cells. Thanks to this second structuring variable, besides compartment dynamics, we aim to represent:
\begin{itemize}
	\item Biochemical reactions inside each compartment. Reaction rates are dependent on local abundances of molecular species, as well as the size (and more generally other physical variable like pH) of the compartment.
	\item Biochemical reactions that occur at the plasma membrane.
	\item Molecular conservation laws between membrane and compartments. Hence we will now distinguish between compartment degradation and compartment recycling.
\end{itemize}
Let $f=f(t,r,a)$ be the population density of compartments at time $t$, size (volume) $r$ and molecular content $a$. Let also $M$ be the molecular quantities at the plasma membrane. The joint evolution equation for $f,M$ is

\begin{multline}\label{eq:coag_frag_mol}\tag{M2a}
\frac{\partial f}{\partial t}+\overbrace{\frac{\partial\left(V(r,a) f(t,r,a)\right)}{\partial a}}^{ \fcolorbox{green}{white}{reactions}}\\
=\underbrace{\frac{1}{2}\int_0^r\int_0^a \kappa(r-r',r')f(t,r-r',a-a')f(t,r',a')\mathrm{d}a'\mathrm{d}r'-\int_0^\infty\int_0^\infty \kappa(r,r')f(t,r,a)f(t,r',a')\mathrm{d}a'\mathrm{d}r'}_{\fcolorbox{red}{white}{coagulation}}\\
+\underbrace{\alpha(r,a,M)}_{\fcolorbox{blue}{white}{endocytosis}}-\underbrace{\gamma(r,a)f(t,r,a)}_{\fcolorbox{blue}{white}{degradation}}-\underbrace{\lambda(r,a)f(t,r,a)}_{\fcolorbox{blue}{white}{recycling}}.
\end{multline}

\begin{equation}\label{eq:coag_frag_mol_membrane}\tag{M2b}
\frac{\mathrm{d}M}{\mathrm{d}t}=\underbrace{J_M(M)}_{ \fcolorbox{green}{white}{reactions}}-\underbrace{\int_0^\infty\int_0^\infty a \alpha(r,a,M)\mathrm{d}a\mathrm{d}r}_{\fcolorbox{blue}{white}{endocytosis}}+\underbrace{\int_0^\infty\int_0^\infty a \lambda(r,a)f(t,r,a)\mathrm{d}a\mathrm{d}r}_{\fcolorbox{blue}{white}{recycling}}.
\end{equation}
In the sequel, we  denote $f(0,r,a) = f_0(r,a)$ and $M(0)=M_0$ the initial conditions of Eqs.~\eqref{eq:coag_frag_mol}-\eqref{eq:coag_frag_mol_membrane}.\\

In Eq.~\eqref{eq:coag_frag_mol}, the extension of the coagulation and endocytosis  processes to include the second structuring variable $a$ is clear. Note that we chose the coagulation kernel to be dependent on the size of the compartments (not their molecular content), for the sake of simplicity. The removal terms are now splitted in two: the degradation occurs at rate $\gamma$, and the recycling occur at rate $\lambda$. The transport term represents the biochemical reactions that modify the molecular content within each compartment, and which occur at rate $V(r,a)$. Also, the endocytosis rate $\alpha$ is necessarily dependent on the molecular content of the plasma membrane $M$, to avoid negative values for $M$, e.g. $\alpha(r,a,0)=0$. %A typical choice (although may be not realistic) could be:
%$$\alpha(x,r,z)=\tilde \alpha(x,r)\indic{x<z}$$
The Eq.~\eqref{eq:coag_frag_mol_membrane} on the scalar variable $M$ was not present in \cite{foretGeneralTheoreticalFramework2012} and, to the best of our knowledge, it is a novelty of our model. 
%seems new.
Still interpreting the molecular content as a quantity of receptors, it allows to represent conservation laws of receptors between plasma membrane and endosomal compartments, giving a satisfactory representation of receptor \textit{trafficking}. In Eq.~\eqref{eq:coag_frag_mol_membrane}, $J_M$ represent biochemical reactions that occur at the plasma membrane, and the two integral terms represent the molecules that are lost of gained at the  plasma membrane through respectively the endocytosis or recycling processes.

In the sequel, the model given by Eqs.~\eqref{eq:coag_frag_mol}-\eqref{eq:coag_frag_mol_membrane} will be referred as our 2D model.

\section{Theoretical properties of the 1D model \eqref{eq:coag_frag}}\label{sec:longtime}

In this section we provide sufficient conditions so that Eq.~\eqref{eq:coag_frag} exhibits a stable steady state. The large-time behaviour of this equation  -also known as \emph{coagulation equation with source and efflux} \cite{chae_existence_1995}, has been studied first in \cite{gajewski_first_1983} with drift, in \cite{dubovskii_mathematical_nodate,laurencot_stationary_2020} without efflux, in \cite{Ofelia2015Thesis} with a bounded coagulation kernel and in  \cite{ghosh_equilibrium_2023} with a singular coagulation kernel. Here, we limit ourselves to give a self-contained proof of exponential stability of the steady state in $L^1$ with bounded kernel. We mainly use a  contraction argument, taking inspiration from \cite{collet_lifshitz-slyozov_1999}.

\begin{thm}\label{theorem_longtime}
  Assume $\alpha$ is integrable, $\kappa$ and $\gamma$ are bounded and morevover $\inf \gamma = \gamma_0 >0$. If
  \begin{equation}\label{hyp:1}
      3 \|\kappa\|_{L^\infty} \|\alpha\|_{L^1} \leq \gamma_0^2\,,
  \end{equation}
  then there exists a unique nonnegative stationary solution in $L^1(\Rb_+)$ of  Eq.~\eqref{eq:coag_frag}, denoted by $f_\infty$. Moreover, 
    \begin{equation*}
\|f_\infty\|_{L^1} \leq \tfrac{\|\alpha\|_{L^1}}{\gamma_0}\,,
\end{equation*}
  and for every solution $f \in C(\Rb_+,L^1(\Rb_+))$ we have
      \begin{equation*}
 \lim_{t\to+\infty} \|f(t) - f_\infty\|_{L^1} = 0\,.
 \end{equation*}
 The convergence is at least exponential with rate $\gamma_0^2-3 \|\kappa\|_{L^\infty} \|\alpha\|_{L^1} >0$.
\end{thm}
We do not expect condition \eqref{hyp:1} to be optimal. See the discussion in Section \ref{sec:4_2}.
\begin{proof}
To simplify the writing, in this section we will denote the coagulation operator by $Q$, namely
\begin{equation*}
    Q(f,f) = \frac{1}{2}\int_0^r \kappa(r-r',r')f(t,r-r')f(t,r')\mathrm{d}r'-\int_0^\infty \kappa(r,r')f(t,r)f(t,r')\mathrm{d}r'\,.
\end{equation*}
 We let $A = \|\alpha\|_{L^1}$, $\kappa_\infty = \|\kappa\|_{L^\infty}$ and $\gamma_\infty = \|\gamma\|_{L^\infty}$.
 
 First we  prove existence and uniqueness of a stationary solution $f_\infty$ thanks to a Banach fixed point argument. Let us define
  \[X = \left\{ f\in L^1(\Rb_+) \, : \, f\geq 0 , \ \|f\|_{L^1} \leq \frac A {\gamma_0}  \right\}. \]
  Consider a constant  $K \geq \gamma_\infty + \kappa_\infty \tfrac A {\gamma_0}$. For $f\in X$ we define
  \[ T_K f =  \frac 1 K \left( \alpha -\gamma f + Q(f,f) + Kf \right).\]
  The operator $T_K$ is well-defined since $Q(f,f)$ is also well-defined for any integrable $f$ provided that $\kappa$ is bounded \cite[Lemma 3]{collet_lifshitz-slyozov_1999}. We aim to apply Banach's fixed point theorem for $T_K$ on $X$ -which is a closed subset of the Banach space $L^1(\Rb_+)$; this will  give a stationary solution of Eq.~\eqref{eq:coag_frag}. Let $f\in X$, we have that 
  \[T_K f \geq  \frac 1 K \left( K - \gamma_\infty   - \kappa_\infty \|f\|_{L^1} \right) f \,. \]
  Indeed, $f$ and $\alpha$ are positive and
  \[Q(f,f) \geq - f(r) \int_0^\infty \kappa(r,r')f(r')\mathrm{d}r'\,.\]
  The fact that $\|f\|_{L^1}\leq \tfrac{A}{\gamma_0}$ and the condition on $K$ entail the positivity of $T_Kf$. Then,
  \begin{equation*}
    \|T_Kf\|_{L^1} = \int_0^\infty \frac 1 K \left( \alpha -\gamma f + Q(f,f) + Kf \right) \leq \frac A K  + \frac{K -\gamma_0}{K}\|f\|\:, 
  \end{equation*}
  noticing that
  \begin{equation*}
     \int_0^\infty Q(f,f)\, \mathrm{d}r  =  - \int_0^\infty \int_0^\infty \kappa(r',r) f(r') f(r) \mathrm{d}r' \mathrm{d}r \leq 0\,.
  \end{equation*}
  But $K \geq \gamma_\infty \geq \gamma_0$, thus
  \[ K\|T_Kf\|_{L^1} \leq A + ( K - \gamma_0) \frac{A}{\gamma_0} = K\frac{A}{\gamma_0}  \]
  and we conclude that $T_Kf$ belongs to $X$. We now show that $T_K$ is a contraction on $X$. Let $f$ and $g$ in $X$. It is straightforward to check that
  \begin{equation*}
 \|T_Kf - T_Kg\|_{L^1} =  (1 -\frac K {\gamma_0}) \|f-g\|_{L^1} + \frac 1 K \|Q(f,f)-Q(g,g)\|_{L^1}\,.
  \end{equation*}
  By a simple computation (see e.g. \cite[Lemma 3]{collet_lifshitz-slyozov_1999} or \cite{Ofelia2015Thesis,Canizo2006Thesis}),

  \begin{equation}
  \label{eq:L1estimate}
  \|Q(f,f)-Q(g,g)\|_{L^1}  \leq \frac{3}{2} \kappa_\infty (\|f\|_{L^1} + \|g\|_{L^1})\|f-g\|_{L^1}\leq  3 \kappa_\infty \frac A {\gamma_0} \|f-g\|_{L^1}.
  \end{equation}
  Thus, we have 
  \[\|T_Kf - T_Kg\|_{L^1} \leq \left( 1 + \frac 1 K ( 3\kappa_\infty \frac A {\gamma_0} -\gamma_0) \right) \|f-g\|_{L^1}\:.\]
  Then the hypothesis Eq.~\eqref{hyp:1} allows us to conclude that there exists a unique fixed point to $T_K$ in $X$.  Moreover, this is the unique stationary solution which is positive and belongs to $L^1$. Indeed, assume we have a positive stationary solution $f$ $\notin X$ i.e. in $L^1$ satisfying $\gamma_0 \int f > A$. Then
  \[0  = \int_0^\infty \left( \alpha -  \gamma f +  Q(f,f)\right) \, dr \leq A - \gamma_0 \int_0^\infty f\, dr <0\:,\]
  which is a contradiction.

  We now turn to the proof of asymptotic stability. Let $f\in C(\Rb_+,L^1(\Rb_+))$ be a nonnegative solution to Eq.~\eqref{eq:coag_frag} in the sense of distributions and let $f_\infty\in L^1(\Rb_+)$ be the stationary solution. We first provide a bound on $f$. We have
  \[ \partial_t \int_0^\infty f(t,r)\mathrm{d}r = A - \gamma_0 \int_0^\infty f(t,r)\mathrm{d}r + \int_0^\infty Q(f,f)(t,r)\mathrm{d}r \leq   A - \gamma_0 \int_0^\infty f(t,r)\mathrm{d}r \,.\]
  Thus,
  \[\|f(t)\|_{L^1}\leq \|f(0)\|_{L^1}e^{-\gamma_0 t} + \frac A {\gamma_0}\,.\]
  It is a classical computation that
  \[ \partial_t |f-f_\infty|  =  - \gamma_0 |f-f_\infty| + (Q(f,f) -Q(f_\infty,f_\infty))\mathop{sign}(f-f_\infty).\]
  Thus, using Eq.~\eqref{eq:L1estimate},
  %by \cite[Lemma 3]{collet_lifshitz-slyozov_1999} 
  \begin{multline*}
  \frac {\mathrm{d}} {\mathrm{d}t} \|f-f_\infty\|_{L^1} \leq \left[ -\gamma_0 +  \frac{3}{2}\kappa_\infty \left(\|f(0)\|_{L^1}e^{-\gamma_0 t} + \frac A {\gamma_0} + \|f_\infty\|_{L^1}\right)\right]\|f-g\|_{L^1} \\
  \leq \left( 3\kappa_\infty \frac A {\gamma_0}  -\gamma_0 +  \frac{3}{2}\kappa_\infty \|f(0)\|_{L^1}e^{-\gamma_0 t}\right)\|f-g\|_{L^1}.
  \end{multline*}
  We conclude that
  \[\|f-f_\infty\|_{L^1} \leq \|f(0)-f_\infty\|_{L^1} e^{\tfrac{2\kappa_\infty}{\gamma_0}\|f(0)\|_{L^1}- \tfrac{\gamma_0^2 - 3\kappa_\infty A}{\gamma_0} t}\,,\]
  which ends the proof.
  \end{proof}

\section{Numerical scheme for 1D and 2D models}\label{sec:numeric}

 In this section, we detail our Finite Volume numerical scheme used to simulate our models, and numerically illustrate their properties. The numerical scheme is adapted from \cite{bourgadeConvergenceFiniteVolume2007,hingantDerivationMathematicalStudy2015}.
 
Again, to simplify the writing, in this section we will denote the coagulation operator by $Q$, namely
\begin{equation*}
    Q(f,f) = \frac{1}{2}\int_0^r\int_0^a \kappa(r-r',r')f(t,r-r',a-a')f(t,r',a')\mathrm{d}a'\mathrm{d}r'-\int_0^\infty\int_0^\infty \kappa(r,r')f(t,r,a)f(t,r',a')\mathrm{d}a'\mathrm{d}r'\,.
\end{equation*}

\subsection{Finite volume scheme}

We want to write a numerical scheme of the 2D model, Eqs.~\eqref{eq:coag_frag_mol}-\eqref{eq:coag_frag_mol_membrane} using a finite volume method. 

First we can remark that for $r>0$ and $a>0$ we can rewrite our Eq.~\eqref{eq:coag_frag_mol} in the following conservative form \cite{bourgadeConvergenceFiniteVolume2007}
\begin{equation}\label{eq:coag_frag_2}
    \begin{split}
        \partial_t f \left(t,r,a\right) = & \dfrac{1}{r a}  \partial_r \partial_a \mathcal{C}\left(f\right)\left(t,r,a\right) - \partial_a \left(V\left(r,a\right) f\left(t,r,a\right)\right) \\ 
        & + \alpha\left(r, a, M\left(t\right)\right) - \lambda\left(r, a\right) f\left(t,r, a\right) - \gamma\left(r, a\right) f\left(t,r,a\right)
    \end{split}
\end{equation}
where 
\begin{equation*}
    \mathcal{C}\left(f\right)\left(t,r,a\right) = \displaystyle \int_0^r \int_0^a r' a' Q \left(f\left(t\right),f\left(t\right)\right)\left(r',a'\right) \mathrm{d}a' \mathrm{d}r'.
\end{equation*}

From Eq.~\eqref{eq:coag_frag_2}, we will first detail the truncation we use (step (i)) and the discretization (step (ii)). Then, we detail the finite volume approximation of the each operator in the right-hand side of Eq.~\eqref{eq:coag_frag_2} (step (iii)). We sum-up the numerical scheme (step (iv)).

\paragraph{\textbf{(i) Truncation}}

As in \cite{hingantDerivationMathematicalStudy2015, bourgadeConvergenceFiniteVolume2007}, to study our equation, we will truncate the size variable to a maximal value $R >0$ and the quantity of reactants variable to a maximal value $A > 0$ and we will choose a truncation of the functional $\mathcal{C}$. We chose the following conservative truncation, given by
\begin{equation*}
    \begin{split}
        \mathcal{C}_{\text{c}}^{RA} \left(f\right)\left(t,r,a\right) = & \frac{1}{2} \displaystyle \int_0^r \int_0^a r' a' \int_0^{r'} \int_0^{a'} \kappa\left(r'-r'',r''\right) f\left(t,r'-r'',a'-a''\right) f\left(t,r'',a''\right) \mathrm{d}a'' \mathrm{d}r'' \mathrm{d}a' \mathrm{d}r' \\ 
        & - \displaystyle \int_0^r \int_0^a r' a' f\left(t,r',a'\right) \int_{0}^{R-r'} \int_{0}^{A-a'} \kappa\left(r',r''\right) f\left(t,r'', a''\right) \mathrm{d}a'' \mathrm{d}r'' \mathrm{d} a'\mathrm{d}r'.
    \end{split}
\end{equation*}

We choose to use a  conservative truncation for our scheme in order to construct a scheme that preserves the conservation properties in the case of  pure coagulation. The equation we will numerically approximate  on the time interval $\left[0,T\right]$, where $T>0$, is thus the following
\begin{equation}\label{eq:coag_frag_simul}
    \begin{split}
        \partial_t f \left(t,r,a\right) = & \dfrac{1}{r a} \partial_r \partial_a \mathcal{C}_{\text{c}}^{RA} \left(f\right)\left(t,r,a\right) - \partial_a \left(V\left(r,a\right) f\left(t,r,a\right)\right) \\
        & + \alpha\left(r, a, M\left(t\right)\right) - \lambda\left(r, a\right) f\left(t,r, a\right) - \gamma\left(r, a\right) f\left(t,r,a\right).
    \end{split}
\end{equation}

\paragraph{\textbf{(ii) Discretization}}
Let $I^r \in \Nb$. We discretize the size interval $\left[0,R\right]$ into $I^r$ intervals. We denote by $\left(r_{i-\frac{1}{2}}\right)_{i \in \left\{1,\dots,I^r+1\right\}}$ a regular mesh of $\left[0,R\right]$ with size step $\Delta r$ and we set
\begin{equation*}
    r_i = \dfrac{r_{i-\frac{1}{2}} + r_{i+\frac{1}{2}}}{2} = \left(i-\dfrac{1}{2}\right) \Delta r, \quad i \in \left\{1,\dots,I^r\right\}.
\end{equation*}

Let $I^a \in \Nb$. We discretize the size interval $\left[0,A\right]$ into $I^a$ intervals. We denote by $\left(a_{j-\frac{1}{2}}\right)_{j \in \left\{1,\dots,I^a+1\right\}}$ a regular mesh of $\left[0,A\right]$ with step $\Delta a$ and we set
\begin{equation*}
    a_j = \dfrac{a_{j-\frac{1}{2}} + a_{j+\frac{1}{2}}}{2} = \left(j-\dfrac{1}{2}\right) \Delta a, \quad j \in \left\{1,\dots,I^a\right\}.
\end{equation*}

For all $i \in \left\{1,\dots,I^r\right\}$ and $j \in \left\{1,\dots,I^a\right\}$, we set 
\begin{equation*}
    \Lambda_{ij} = \left[ r_{i-\frac{1}{2}}, r_{i+\frac{1}{2}}\right] \times \left[ a_{j-\frac{1}{2}}, a_{j+\frac{1}{2}}\right].
\end{equation*}

 Let $\Delta t >0$ be the time step. We discretize $\left[0,T\right]$ by the set of points $\left\{ t^n = n \, \Delta t,\; n \in \left\{0, \dots,N \right\} \right\}$, where $N =  \left\lfloor\frac{T}{\Delta t}\right\rfloor$.

\paragraph{\textbf{(iii) Finite volume approximation}}

For all $n \in \left\{0, \dots, N \right\}$, $i \in \left\{1, \dots, I^r\right\}$ and $j \in \left\{1,\dots,I^a\right\}$, we denote $f_{i,j}^n$ an approximation of the function $f$ at the  point $\left(t^n, r_i,a_j\right)$. We set
\begin{equation*}
    f_{i,j}^n = \dfrac{1}{\Delta r \Delta a} \int_{\Lambda_{ij}} f\left(t^n, r, a\right) \mathrm{d}a \mathrm{d}r.
\end{equation*}

We denote by $M^n$ an approximation of $M\left(t^n\right)$.
First, we write the explicit time Euler scheme associated to the Eq.~\eqref{eq:coag_frag_simul}. For all $n \in \left\{0, \dots,N-1 \right\}$, we have
\begin{equation}\label{eq:coag_frag_euler}
    \begin{split}
        \dfrac{f \left(t^{n+1},r,a\right) - f \left(t^n,r,a\right)}{\Delta t} = & \dfrac{1}{r a}  \partial_r \partial_a \mathcal{C}_{\text{c}}^{RA}\left(f\right)\left(t^n,r,a\right) - \partial_a \left(V\left(r,a\right) f\left(t^n,r,a\right)\right) \\ 
        & + \alpha\left(r, a, M^n\right) - \lambda\left(r, a\right) f\left(t^n,r, a\right) - \gamma\left(r, a\right) f\left(t^n,r,a\right).
    \end{split}
\end{equation}

Then, for all $i \in \left\{1,\dots,I^r\right\}$ and $j \in \left\{1,\dots,I^a\right\}$, integrating Eq.~\eqref{eq:coag_frag_euler} over $\Lambda_{ij}$ and using some approximations, we have:
\begin{equation*}\label{eq:coag_frag_euler_int}
    \begin{split}
        \dfrac{\Delta r \Delta a}{\Delta t} \left(f_{i,j}^{n+1}-f_{i,j}^n\right) = & \dfrac{1}{r_{i-\frac{1}{2}} a_{j-\frac{1}{2}}} \left( \mathcal{C}_{i,j}^n - \mathcal{C}_{i-1,j}^n - \mathcal{C}_{i,j-1}^n + \mathcal{C}_{i-1,j-1}^n \right) \\
        & - \Delta r \cdot \left[ A^{\text{up}} \left(W_{i,j+1}, f_{i,j}^n, f_{i,j+1}^n\right) - A^{\text{up}} \left(W_{i,j}, f_{i,j-1}^n, f_{i,j}^n\right) \right] \\ 
        & + \int_{\Lambda_{ij}} \alpha\left(r, a, M^n\right) \mathrm{d}a \mathrm{d}r - f_{i,j}^n  \int_{\Lambda_{ij}} \lambda\left(r, a\right) \mathrm{d}a \mathrm{d}r - f_{i,j}^n  \int_{\Lambda_{ij}} \gamma\left(r, a\right) \mathrm{d}a \mathrm{d}r
    \end{split}
\end{equation*}
where $\mathcal{C}_{i,j}^n$ is an approximation of $\mathcal{C}_{\text{c}}^{RA}\left(f\right)\left(t^n,r_{i+\frac{1}{2}},a_{j+\frac{1}{2}}\right)$, $W_{i,j} = \dfrac{1}{2}\left(V\left(r_{i+\frac{1}{2}}, a_{j-\frac{1}{2}}\right) + V\left(r_{i-\frac{1}{2}}, a_{j-\frac{1}{2}}\right)\right)$ and $f_{i,0}^n = f_{i,I^r+1}^n = 0$.
We define the operator $A^\text{up}$, used for an upwind approximation of the transport term, as follows
\begin{equation*}
    A^{\text{up}}\left(u, f_+, f_-\right) = \left\{
    \begin{array}{ll}
        u f_+ & \text{if } u\geq0 \\
        u f_- & \text{if } u<0.
    \end{array}
    \right.
\end{equation*}

Long but straightforward  calculations show that
\begin{small}
\begin{equation*}
    \begin{split}
        \mathcal{C}_{\text{c}}^{RA}\left(f\right)\left(t^n,r_{i+\frac{1}{2}},a_{j+\frac{1}{2}}\right) & = \frac{1}{2}  \int_0^{r_{i+\frac{1}{2}}} \int_0^{a_{j+\frac{1}{2}}} \int_0^{r'} \int_0^{a'} r' a' \kappa\left(r'-r'',r''\right)  f\left(t^n,r'-r'',a'-a''\right) f\left(t^n,r'',a''\right) \mathrm{d}a'' \mathrm{d}r'' \mathrm{d}a' \mathrm{d}r' \\
        & \quad - \int_0^{r_{i+\frac{1}{2}}} \int_0^{a_{j+\frac{1}{2}}} r' a' f\left(t^n,r',a'\right) \int_{0}^{R-r'} \int_{0}^{A-a'} \kappa\left(r',r''\right) f\left(t^n,r'', a''\right) \mathrm{d}a'' \mathrm{d}r'' \mathrm{d} a'\mathrm{d}r' \\
        & = \frac{1}{2} \displaystyle \sum_{k=1}^i \sum_{m=1}^j \int_{r_{k-\frac{1}{2}}}^{r_{k+\frac{1}{2}}} \int_{a_{m-\frac{1}{2}}}^{a_{m+\frac{1}{2}}} \int_0^{r_{i+\frac{1}{2}}-r'} \int_0^{a_{j+\frac{1}{2}}-a'}\left(a'+a''\right)\left(r'+r''\right) \kappa\left(r',r''\right) f\left(t^n,r',a'\right) \\
        & \quad \quad f\left(t^n,r'',a''\right) \mathrm{d}a'' \mathrm{d}r'' \mathrm{d}a' \mathrm{d}r' \\
        & \quad - \displaystyle \sum_{k=1}^i \sum_{m=1}^j \int_{r_{k-\frac{1}{2}}}^{r_{k+\frac{1}{2}}} \int_{a_{m-\frac{1}{2}}}^{a_{m+\frac{1}{2}}} r' a' f\left(t^n,r',a'\right) \int_{0}^{R-r'} \int_{0}^{A-a'} \kappa\left(r',r''\right) f\left(t^n,r'', a''\right) \mathrm{d}a'' \mathrm{d}r'' \mathrm{d} a'\mathrm{d}r' \\
        & \simeq \frac{1}{2} \left(\Delta r \Delta a\right)^2 \displaystyle \sum_{k=1}^i \sum_{m=1}^j \sum_{k'=1}^{i-k+1} \sum_{m'=1}^{j-m+1} \left(a_{m-\frac{1}{2}}+a_{m'-\frac{1}{2}}\right)\left(r_{k-\frac{1}{2}}+r_{k'-\frac{1}{2}}\right) \kappa_{k,k'} f_{k,m}^n f_{k',m'}^n  \\
        & \quad - \left(\Delta r \Delta a\right)^2\displaystyle \sum_{k=1}^i \sum_{m=1}^j \sum_{k'=1}^{I^r-k+1} \sum_{m'=1}^{I^a-m+1} a_{m-\frac{1}{2}}r_{k-\frac{1}{2}} \kappa_{k,k'} f_{k,m}^n f_{k',m'}^n
    \end{split}
\end{equation*}
\end{small}
where $\kappa_{k,k'}$ is an approximation of $\kappa\left(r,r'\right)$ with $r\in \left[r_{k-\frac{1}{2}}, r_{k+\frac{1}{2}}\right]$ and $r'\in \left[r_{k'-\frac{1}{2}}, r_{k'+\frac{1}{2}}\right]$.

Then we set
\begin{equation*}
    \begin{split}
        \mathcal{C}_{i,j}^n & = \frac{1}{2} \left(\Delta r \Delta a\right)^2 \displaystyle \sum_{k=1}^i \sum_{m=1}^j \sum_{k'=1}^{i-k+1} \sum_{m'=1}^{j-m+1} \left(a_{m-\frac{1}{2}}+a_{m'-\frac{1}{2}}\right)\left(r_{k-\frac{1}{2}}+r_{k'-\frac{1}{2}}\right) \kappa_{k,k'} f_{k,m}^n f_{k',m'}^n  \\ 
        & \quad - \left(\Delta r \Delta a\right)^2\displaystyle \sum_{k=1}^i \sum_{m=1}^j \sum_{k'=1}^{I^r-k+1} \sum_{m'=1}^{I^a-m+1} a_{m-\frac{1}{2}}r_{k-\frac{1}{2}} \kappa_{k,k'} f_{k,m}^n f_{k',m'}^n
    \end{split}
\end{equation*}
and we have
\begin{equation*}
    \begin{split}
        \mathcal{C}_{i,j}^n - \mathcal{C}_{i-1,j}^n - \mathcal{C}_{i,j-1}^n + \mathcal{C}_{i-1,j-1}^n & = \frac{1}{2} \left(\Delta r \Delta a\right)^2 a_{j-\frac{1}{2}} r_{i-\frac{1}{2}} \displaystyle \sum_{k=1}^i \sum_{m=1}^j  \kappa_{k,i-k+1} f_{k,m}^n f_{i-k+1,j-m+1}^n \\
        & \quad - \left(\Delta r \Delta a\right)^2 a_{j-\frac{1}{2}}r_{i-\frac{1}{2}} f_{i,j}^n \displaystyle \sum_{k=1}^{I^r-i+1} \sum_{m=1}^{I^a-j+1} \kappa_{i,k} f_{k,m}^n.
    \end{split}
\end{equation*}

It remains to find a way to compute $M^n$. We use an explicit time Euler scheme to do that. Then we have
\begin{equation*}
    \begin{split}
        \dfrac{M^{n+1}-M^n}{\Delta t} = & J_M\left(M^n\right) - \int_0^R \int_0^A a'\alpha\left(r',a',M^n\right) \mathrm{d}a' \mathrm{d}r' \\
        & + \sum_{i=1}^{I^r} \sum_{j=1}^{I^a} \int_{r_{i-\frac{1}{2}}}^{r_{i+\frac{1}{2}}} \int_{a_{j-\frac{1}{2}}}^{a_{j+\frac{1}{2}}} a' \lambda\left(r',a',M^n\right) f_{i,j}^n \mathrm{d}a' \mathrm{d}r'.
    \end{split}
\end{equation*}

\paragraph{\textbf{(iv) Finite volume scheme}}
Finally, the scheme of the model \eqref{eq:coag_frag_mol}-\eqref{eq:coag_frag_mol_membrane} is given by:
\begin{itemize}
    \item Initialization
    for $i \in \left\{1,\dots,I^r\right\}$ and $j \in \left\{1,\dots,I^a\right\}$ we set
    \begin{equation*}
        f_{i,j}^0 = \dfrac{1}{\Delta r \Delta a} \int_{\Lambda_{ij}} f\left(0, r, a\right) \mathrm{d}a \mathrm{d}r.
    \end{equation*}
    
    \item Time iteration
    for all $n \in \left\{ 1, \dots, N\right\}$, $i \in \left\{1,\dots,I^r\right\}$ and $j \in \left\{1,\dots,I^a\right\}$, we have
    \begin{equation}\label{eq:coag_frag_scheme}
        \begin{split}
            f_{i,j}^{n+1} = & f_{i,j}^n + \Delta t \Delta r \Delta a \left[ \frac{1}{2} \displaystyle \sum_{k=1}^i \sum_{m=1}^j \kappa_{k,i-k+1} f_{k,m}^n f_{i-k+1,j-m+1}^n -  \displaystyle f_{i,j}^n \sum_{k=1}^{I^r-i+1} \sum_{m=1}^{I^a-j+1} \kappa_{i,k} f_{k,m}^n \right] \\ 
            & - \dfrac{\Delta t}{\Delta a} \left[ A^{\text{up}} \left(W_{i,j+1}, f_{i,j}^n, f_{i,j+1}^n\right) - A^{\text{up}} \left(W_{i,j}, f_{i,j-1}^n, f_{i,j}^n\right) \right] \\
            & + \dfrac{\Delta t}{\Delta r \Delta a} \left[ \int_{\Lambda_{ij}} \alpha\left(r, a, M^n\right) \mathrm{d}a \mathrm{d}r -   f_{i,j}^n  \int_{\Lambda_{ij}} \lambda\left(r, a\right) \mathrm{d}a \mathrm{d}r - f_{i,j}^n  \int_{\Lambda_{ij}} \gamma\left(r, a\right) \mathrm{d}a \mathrm{d}r\right].\\
        M^{n+1} = & M^n + \Delta t \, J_M\left(M^n\right) - \Delta t \int_0^R \int_0^A a'\alpha\left(r',a',M^n\right) \mathrm{d}a' \mathrm{d}r' \\
        & + \Delta t \sum_{i=1}^{I^r} \sum_{j=1}^{I^a} \int_{r_{i-\frac{1}{2}}}^{r_{i+\frac{1}{2}}} \int_{a_{j-\frac{1}{2}}}^{a_{j+\frac{1}{2}}} a' \lambda\left(r',a',M^n\right) f_{i,j}^n \mathrm{d}a' \mathrm{d}r'.
    \end{split}
\end{equation}
\end{itemize}

\begin{rmk}
    With the same tools, we can write a scheme in dimension 1 noticing that if we set $g\left(t,r\right) = r f\left(t,r\right)$ we have
    \begin{equation*}
        \partial_t g\left(t, r\right) = -\partial_r J\left(f\right)\left(t,r\right) + r \alpha\left(r,M\left(t\right)\right) - \lambda\left(r\right) g\left(t,r\right) - \gamma\left(r\right) g\left(t,r\right)
    \end{equation*}
    where $J\left(f\right)$ is defined as follows
    \begin{equation*}
        J\left(f\right)\left(t,r\right) = \displaystyle \int_0^r \int_{r-r'}^\infty r' \kappa\left(r',r''\right) f\left(t,r'\right) f\left(t,r''\right) \mathrm{d}r'' \mathrm{d}r'.
    \end{equation*}
    A conservative truncation is then given by the following:
    \begin{equation*}
        J_c^R \left(f\right)\left(t,r\right) = \displaystyle \int_0^r \int_{r-r'}^{R-r'} r' \kappa\left(r',r''\right) f\left(t,r'\right) f\left(t,r''\right) \mathrm{d}r'' \, \mathrm{d}r'.
    \end{equation*}
\end{rmk}

\subsection{Numerical checking - Long-time behavior/stability}
\label{sec:4_2}    
    In this section, we will investigate the behavior of the numerical scheme given by Eq.~\eqref{eq:coag_frag_scheme} in some particular cases.
    
    Our criteria will be based on some moments evaluation that we introduce now.  First, we will evaluate the numerical scheme in the case of a pure coagulation model ($\alpha=\gamma=\lambda=V=J_M=0$). We consider the pure coagulation equation:
    \begin{equation} \label{eq:coag}
        \partial_t f\left(t,r,a\right) = Q\left(f\left(t\right),f\left(t\right)\right) \left(r,a\right), \quad \text{with } t>0, r>0 \text{ and } a>0.
    \end{equation}
 Let $Q_{i,j}^n\left(f,f\right)$ be the discrete coagulation operator,
    \begin{equation*}
        Q_{i,j}^n\left(f,f\right) = \Delta r \Delta a \left(\frac{1}{2} \displaystyle \sum_{k=1}^i \sum_{m=1}^j  \kappa_{k,i-k+1} f_{k,m}^n f_{i-k+1,j-m+1}^n - f_{i,j}^n \displaystyle \sum_{k=1}^{I^r-i+1} \sum_{m=1}^{I^a-j+1} \kappa_{i,k} f_{k,m}^n\right).
    \end{equation*}
    The discrete equation associated to Eq.~\eqref{eq:coag} is given by
    \begin{equation} \label{eq:coag_dis}
        \dfrac{f_{i,j}^{n+1}-f_{i,j}^n}{\Delta t} = Q_{i,j}^n\left(f,f\right), \quad n \in \left\{ 1, \dots, N\right\}, i \in \left\{1,\dots,I^r\right\} \text{ and } j \in \left\{1,\dots,I^a\right\}.
    \end{equation}
    We remark that we keep the same properties on the moments dynamics for the two Eqs.~\eqref{eq:coag} and \eqref{eq:coag_dis}. For some test function $\varphi$, we define $H\left(\varphi, t\right)$ the moment of $f$ associated with the function $\varphi$ at time $t>0$
    \[H\left(\varphi, t\right) := \int_0^\infty \int_0^\infty \varphi\left(r,a\right) f\left(t,r,a\right) \mathrm{d}a \mathrm{d}r\,, \]
    and $H^n\left(\varphi\right)$ its discrete analogue moment, associated with the function $\varphi$ at time $t^n$ with $n \in \left\{1,\dots,N\right\}$. Then
    \begin{equation*}
        \begin{split}
            \dfrac{\mathrm{d}}{\mathrm{d}t} H\left(\varphi, t\right) & = \dfrac{\mathrm{d}}{\mathrm{d}t} \left[ \int_0^\infty \int_0^\infty \varphi\left(r,a\right) f\left(t,r,a\right) \mathrm{d}a \mathrm{d}r \right] \\
            & = \int_0^\infty \int_0^\infty \varphi\left(r,a\right) Q\left(f\left(t\right),f\left(t\right)\right) \left(r,a\right) \mathrm{d}a \mathrm{d}r \\
            & = \frac{1}{2} \int_0^\infty \int_0^\infty \int_0^\infty \int_0^\infty \left[\varphi\left(r+r',a+a'\right) - \varphi\left(r,a\right) - \varphi\left(r',a'\right) \right] \kappa\left(r,r'\right) f\left(t,r,a\right) f\left(t,r',a'\right) \mathrm{d}a' \mathrm{d}r' \mathrm{d}a \mathrm{d}r
        \end{split}
    \end{equation*}
    and we have
    \begin{equation*}
        \begin{split}
            \dfrac{H^{n+1}\left(\varphi\right)-H^n\left(\varphi\right)}{\Delta t} & = \Delta r \Delta a \sum_{i=1}^{I^r} \sum_{j=1}^{I^a} \varphi\left(r_{i-\frac{1}{2}},a_{j-\frac{1}{2}}\right) \dfrac{f_{i,j}^{n+1}-f_{i,j}^n}{\Delta t} \\
            & = \Delta r \Delta a \sum_{i=1}^{I^r} \sum_{j=1}^{I^a} \varphi\left(r_{i-\frac{1}{2}},a_{j-\frac{1}{2}}\right) Q_{i,j}^n\left(f,f\right) \\
            & = \frac{\left(\Delta r \Delta a\right)^2}{2} \sum_{i=1}^{I^r} \sum_{j=1}^{I^a} \sum_{k=1}^{I^r-i+1} \sum_{m=1}^{I^a-j+1} \Big{[} \varphi\left(r_{i+k-1-\frac{1}{2}},a_{j+m-1-\frac{1}{2}}\right) \\ & \quad - \varphi\left(r_{i-\frac{1}{2}},a_{j-\frac{1}{2}}\right) - \varphi\left(r_{k-\frac{1}{2}},a_{m-\frac{1}{2}}\right) \Big{]} \kappa_{i,k} f_{i,j}^n f_{k,m}^n.
        \end{split}
    \end{equation*}
    
    Thus choosing $\varphi\left(r,a\right) = r$ or $\varphi\left(r,a\right) = a$, we have that the first-order moments are constant in time both at the discrete and continuous levels. We also have that the zeroth-order moments are non-increasing functions of time in both cases. For the pure coagulation model given by Eq.~\eqref{eq:coag}, choosing an affine kernel $\kappa\left(r,r'\right) = K_0 + K_1\left(r+r'\right)$ leads to a closed moment equation in the form of an ODE system. Indeed, in such case, it is easy to see that we have the following ODE system for the moments of order 0 and 1: %\textcolor{teal}{and 2}:
    \begin{equation}\label{eq:coag_moments_ODE}
        \left\{
        \begin{array}{l}
            \dfrac{\mathrm{d}}{\mathrm{d}t} H\left(1,t\right) = - \dfrac{1}{2}K_0 \left(H\left(1,t\right)\right)^2 - K_1 H\left(1,t\right) H\left(r,t\right), \\ \\
            \dfrac{\mathrm{d}}{\mathrm{d}t} H\left(r,t\right) = 0, \\ \\
            \dfrac{\mathrm{d}}{\mathrm{d}t} H\left(a,t\right) = 0. \\
            %\dfrac{\mathrm{d}}{\mathrm{d}t} H\left(r^2,t\right) = K_0 \left(H\left(r,t\right)\right)^2 +2 K_1 H\left(r,t\right) H\left(r^2,t\right) \\
            %\dfrac{\mathrm{d}}{\mathrm{d}t} H\left(a^2,t\right) = K_0 \left(H\left(a,t\right)\right)^2 +2 K_1 H\left(a,t\right) H\left(ar,t\right) \\
            %\dfrac{\mathrm{d}}{\mathrm{d}t} H\left(ar,t\right) = K_0 H\left(a,t\right)H\left(r,t\right) + K_1 H\left(r,t\right) H\left(ar,t\right) + K_1 H\left(a,t\right) H\left(r^2,t\right)
        \end{array}
        \right.
    \end{equation}

    Numerical solutions of Eq.~\eqref{eq:coag_moments_ODE} are computed using the standard ODE solver ODEProblem from the Julia package DifferentialEquations, and are compared with moments calculated from our numerical scheme \eqref{eq:coag_frag_scheme} of the 2D model, computed as
    \[ H^n(h) = \Delta a \Delta r \sum_{i=1}^{I^r} \sum_{j=1}^{I^a}  h\left(r_{i-\frac{1}{2}},a_{j-\frac{1}{2}}\right) f_{i,j}^n. \]
    
\begin{figure}[h!]
	\includegraphics[width=0.45\linewidth]{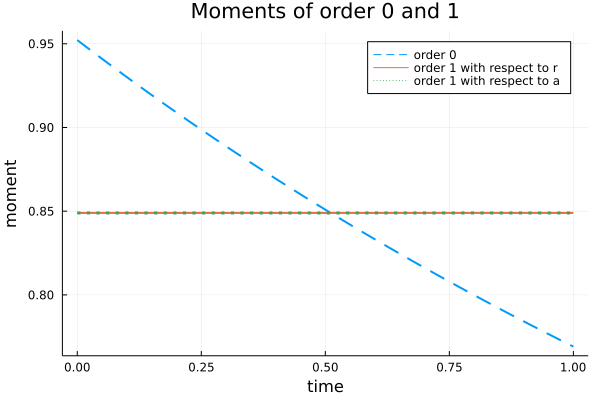}
	\includegraphics[width=0.45\linewidth]{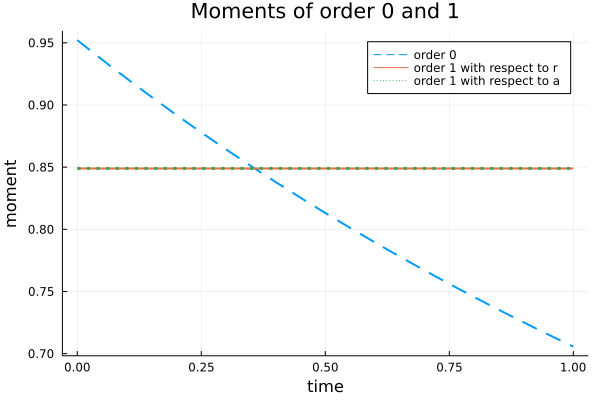}
	\vspace{-0.5cm}
	\caption{{\small\textit{Moments of order $0$ and $1$ computed from our numerical scheme \eqref{eq:coag_frag_scheme}, in the pure coagulation case for two different kernels $\kappa \left(r,r'\right) = 0.5$ for the left picture and $\kappa \left(r,r'\right) = 0.5 + 0.1 \left(r+r'\right)$ for the right one. The results are obtained with the following parameters: $\Delta t = 10^{-4}$, $I^r = I^a = 40$, $R=A=10$ and $f_0\left(r,a\right) =0.5\cdot \left[\mathcal{N}_r(1.5,0.15)\times\mathcal{N}_a(0.5,0.3)+\mathcal{N}_r(0.5,0.3)\times \mathcal{N}_a(1.5,0.15)\right]$. { Here} $\mathcal{N}(\mu,\sigma)$ {stands for a Gaussian density function of mean $\mu$ and standard deviation $\sigma$.}  }}}\label{fig:moments}
	\vspace{-0.3cm}
    \end{figure}

\begin{figure}[h!]
	\includegraphics[width=0.45\linewidth]{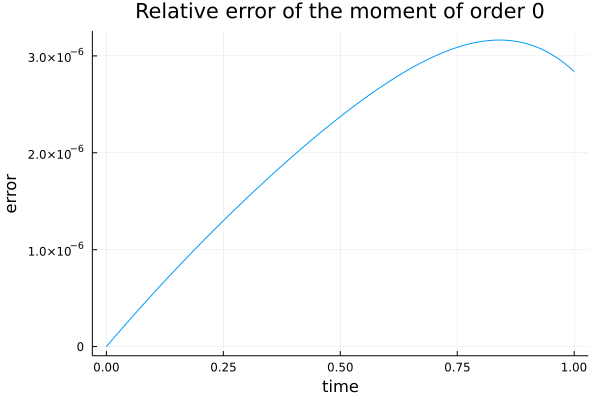}
	\includegraphics[width=0.45\linewidth]{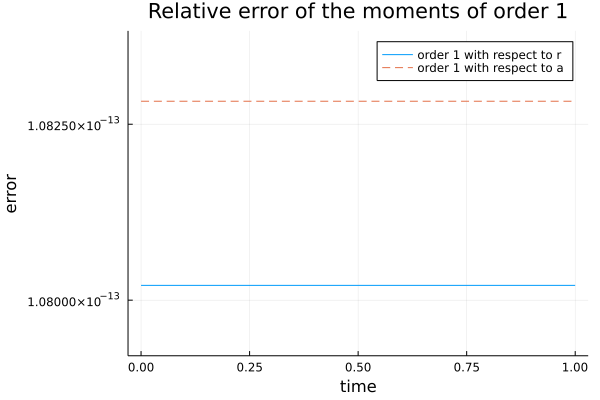} 
	\vspace{-0.5cm}
	\caption{{\small\textit{Relative error of the moments of order 0 and 1 between moments computed from our numerical scheme \eqref{eq:coag_frag_scheme} and from the ODE system Eq.~\eqref{eq:coag_moments_ODE}, in the case of a constant kernel $\kappa \left(r,r'\right) = 0.5$ with the same parameters as in Figure \ref{fig:moments}. }}}\label{fig:error_moments_constant}
	\vspace{-0.3cm}
    \end{figure}

\begin{figure}[h!]
	\includegraphics[width=0.45\linewidth]{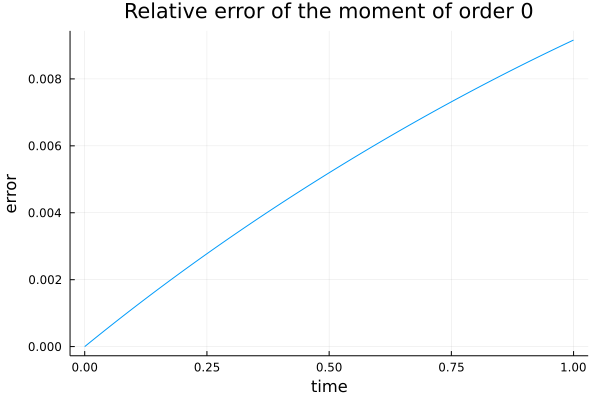}
	\includegraphics[width=0.45\linewidth]{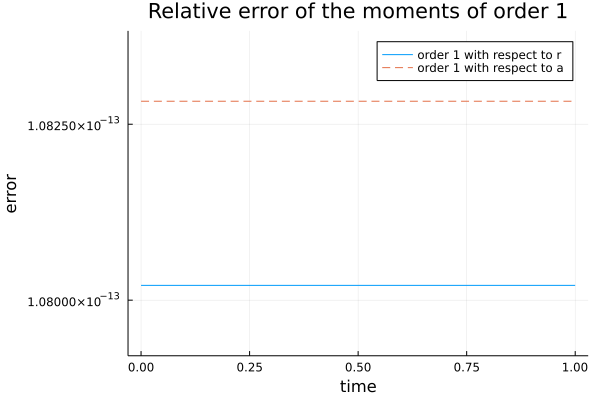} 
	\vspace{-0.5cm}
	\caption{{\small\textit{Relative error of the moments of order 0 and 1 between moments computed from our numerical scheme \eqref{eq:coag_frag_scheme} and from the ODE system Eq.~\eqref{eq:coag_moments_ODE}, in the case of an affine kernel $\kappa \left(r,r'\right) = 0.5 + 0.1 \left(r+r'\right)$ with the same parameters as in figure \ref{fig:moments}. }}}\label{fig:error_moments_linear}
	\vspace{-0.3cm}
    \end{figure}

     As we expect from Eq.~\eqref{eq:coag_moments_ODE}, we recover from our numerical scheme \eqref{eq:coag_frag_scheme} that the moment of order 0 is a nonincreasing function and the moments of order 1 are constant (Figure \ref{fig:moments}). Moreover, we also observed that these numerical results are very close to the numerical solutions directly computed from Eq.~\eqref{eq:coag_moments_ODE} (Figures \ref{fig:error_moments_constant} and \ref{fig:error_moments_linear}). Relative error of moments of order $0$ are increasing through time, as expected from the fact that the size-truncation of the numerical scheme \eqref{eq:coag_frag_scheme} implies more and more error as compartments gets bigger in pure-coagulation dynamics.
    
\begin{figure}[h!]
	\includegraphics[width=0.3\linewidth]{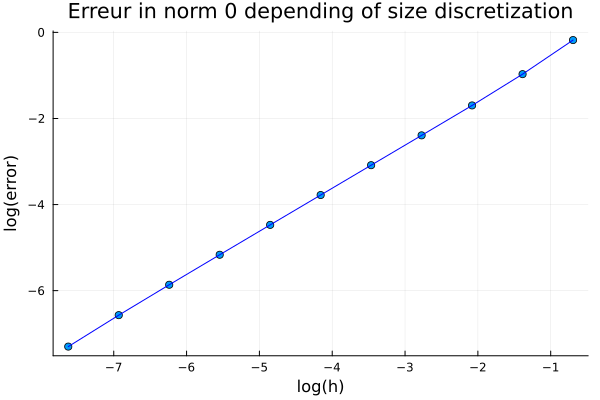}
	\includegraphics[width=0.3\linewidth]{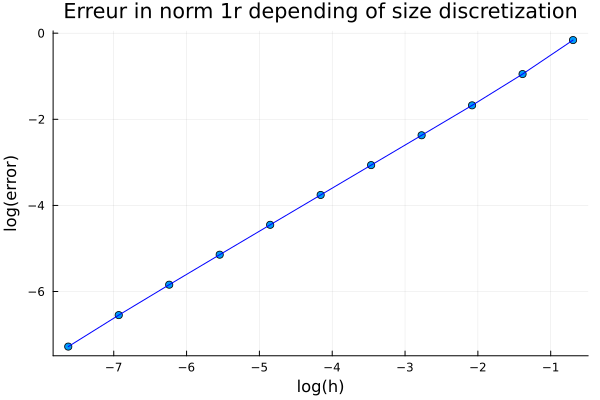} 
	\includegraphics[width=0.3\linewidth]{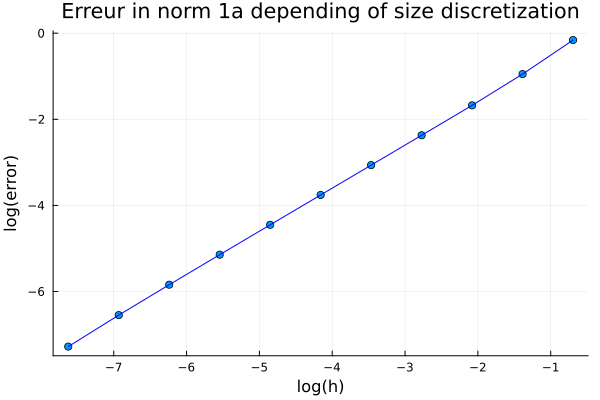} 
	\includegraphics[width=0.3\linewidth]{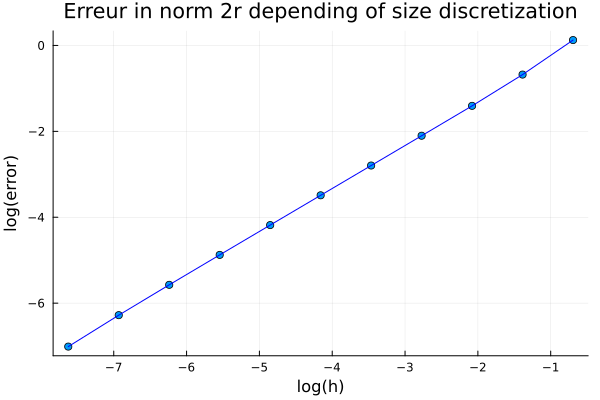} 
	\includegraphics[width=0.3\linewidth]{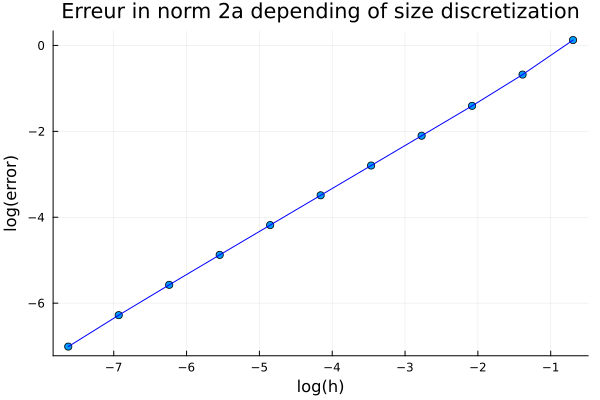} 
	\includegraphics[width=0.3\linewidth]{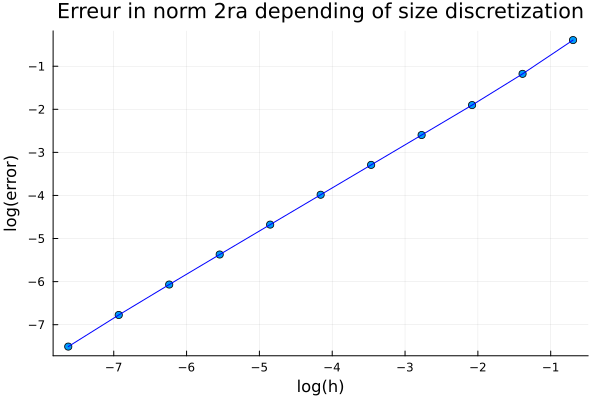}
	\vspace{-0.2cm}
	\caption{{\small\textit{Illustration of the convergence of the scheme in the case of the pure coagulation. The results are obtained with the following parameters: $\Delta t = 0.0005$, $I^r = I^a = 6$, $R=A=3$, $\kappa\left(r,r'\right) = 0.5$ and $f_0\left(r,a\right) = 0.5 \cdot \left[\mathcal{N}_r\left(1.5,0.15\right)\times \mathcal{N}_a\left(0.5,0.3\right) + \mathcal{N}_r\left(0.5,0.3\right)\times \mathcal{N}_a\left(1.5,0.15\right) \right]$. }}}\label{fig:conv_coag_pure}
	\vspace{-0.3cm}
    \end{figure}

\begin{figure}[h!]
	\includegraphics[width=0.3\linewidth]{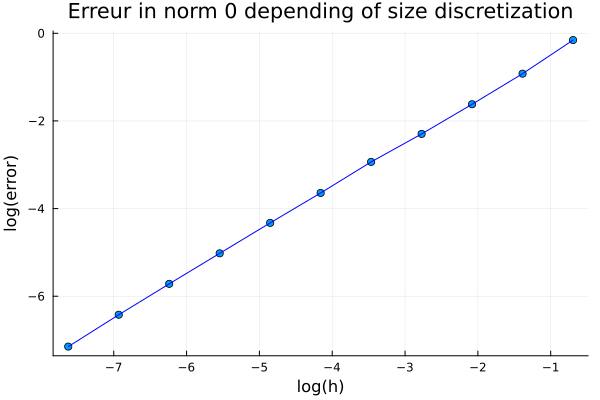}
	\includegraphics[width=0.3\linewidth]{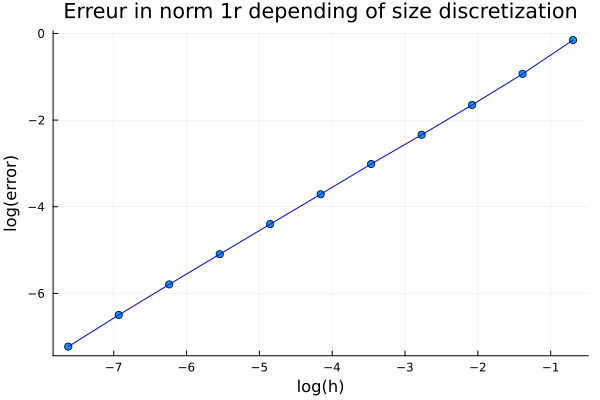} 
	\includegraphics[width=0.3\linewidth]{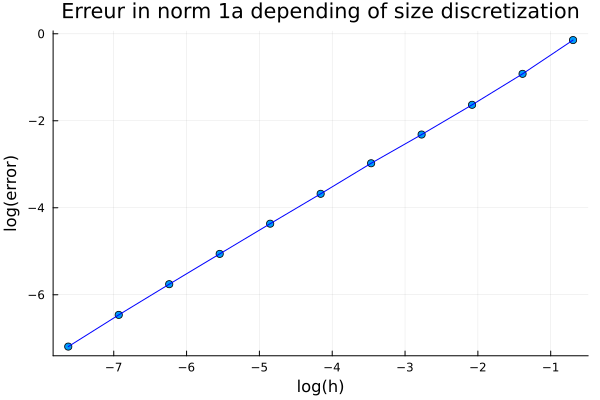} 
	\includegraphics[width=0.3\linewidth]{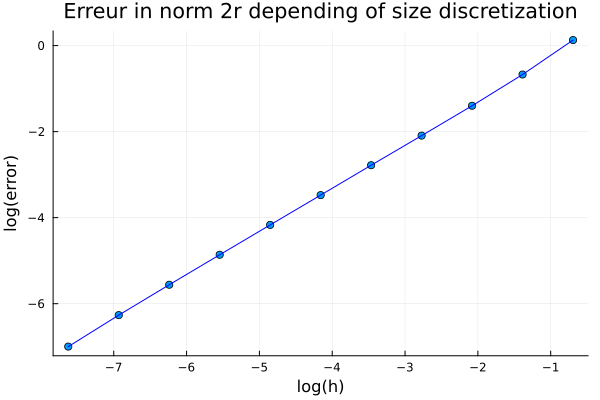} 
	\includegraphics[width=0.3\linewidth]{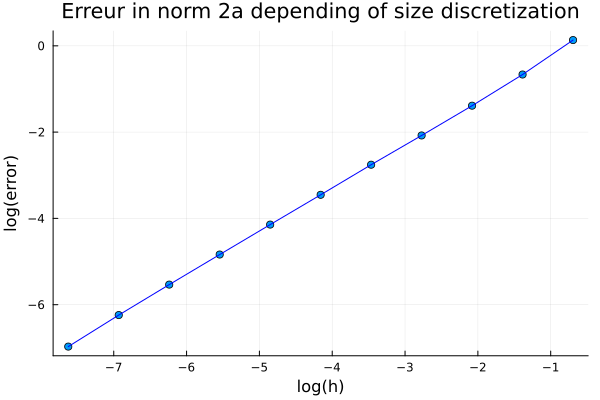} 
	\includegraphics[width=0.3\linewidth]{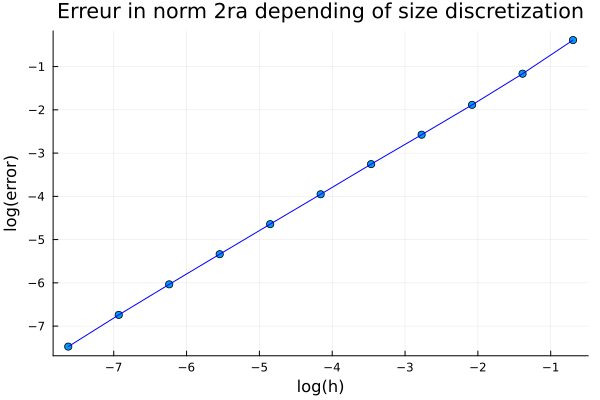}
	\vspace{-0.2cm}
	\caption{{\small\textit{Illustration of the convergence of the scheme in the case of the general case. The results are obtained with the following parameters: $\Delta t = 0.0005$, $I^r = I^a = 6$, $R=A=3$, $\kappa\left(r,r'\right) = 0.5$, $f_0\left(r,a\right) = 0.5 \cdot \left[\mathcal{N}_r\left(1.5,0.15\right)\times \mathcal{N}_a\left(0.5,0.3\right) + \mathcal{N}_r\left(0.5,0.3\right)\times \mathcal{N}_a\left(1.5,0.15\right) \right]$, $\alpha \left(r,a,M\right) = 0.1 \cdot M \cdot \left[\mathcal{N}_r\left(0.6,0.01\right)\times \mathcal{N}_a\left(0.3,0.05\right) + \mathcal{N}_r\left(0.3,0.05\right)\times \mathcal{N}_a\left(0.6,0.01\right) \right]$, $\gamma \left(r,a\right) = 20 \left(r-5\right)^4 \indic{r>5} + 10^{-5} $,  $\lambda\left(r,a\right) = 10^{-2}\cdot \mathcal{P}_r(10,0)$, $V\left(r,a\right) = 0$, $J_M\left(M\right) = 0$ and $M\left(0\right) = 20$. }}}\label{fig:conv_coag_reac}
	\vspace{-0.3cm}
    \end{figure}

    {We now investigate the convergence property of our numerical scheme given by Eq.~\eqref{eq:coag_frag_scheme} in a more general setting. To do this, we fix a small time step $\Delta t$ and the discretization step $\Delta r$ and $\Delta a$. We assume, in this case, $\Delta r$ and $\Delta a$ to be equal and we denote by $h$ the discretization parameter. We have $h=\Delta r = \Delta a$. The finite volume is then of size $h^2$. We compute the solution at time $\Delta t$. Then we divide $h$ by two and compute the solution at time $\Delta t$ with the finite volume size $\left(\frac{h}{2}\right)^2$. We repeat it until reaching the size $\left(\frac{h}{2^{11}}\right)^2$.} We choose this finest solution as reference solution and compare the solutions defined on the  coarser grids to this reference solution with the following different norms, defined for a function $g$ as follows:
    \begin{align*}
        \left\|g\right\|_0  & = \int_0^R \int_0^A \left|g\left(r,a\right)\right| \mathrm{d}a \mathrm{d}r, \\
        \left\|g\right\|_{1r}  & = \int_0^R \int_0^A r \cdot \left|g\left(r,a\right)\right| \mathrm{d}a \mathrm{d}r, \\
        \left\|g\right\|_{1a}  & = \int_0^R \int_0^A a \cdot \left|g\left(r,a\right)\right| \mathrm{d}a \mathrm{d}r, \\
        \left\|g\right\|_{2r}  & = \int_0^R \int_0^A r^2 \cdot \left|g\left(r,a\right)\right| \mathrm{d}a \mathrm{d}r, \\
        \left\|g\right\|_{2a}  & = \int_0^R \int_0^A a^2 \cdot \left|g\left(r,a\right)\right| \mathrm{d}a \mathrm{d}r, \\
        \left\|g\right\|_{2ra}  & = \int_0^R \int_0^A r a \cdot \left|g\left(r,a\right)\right| \mathrm{d}a \mathrm{d}r.
    \end{align*}
    We compute this in the pure coagulation case Eq.~\eqref{eq:coag} (Figure \ref{fig:conv_coag_pure}) and on the general case Eqs.~\eqref{eq:coag_frag_mol}-\eqref{eq:coag_frag_mol_membrane} (Figure \ref{fig:conv_coag_reac}), with reactions, endocytosis, degradation and recycling as detailed in the legend of Figure \ref{fig:conv_coag_reac}, where  we introduce the notations $\mathcal{N}_x(\mu,\sigma):x \mapsto\frac{1}{\sqrt{2 \pi }\sigma}e^{-\frac{1}{2}\left( \frac{x-\mu}{\sigma}\right)^2}$  and $\mathcal{P}_x(\bar{X},\epsilon):x \mapsto \left(\frac{\bar{X}+\epsilon-x}{\bar{X}+\epsilon}\right)^{\frac{1}{3}} \left(\frac{x}{\bar{X}}\right)^{\frac{2}{3}}$.\\
      
    In both cases, we observe that the error decreases linearly (in log-log scale) when the size of the discretization step $h$ does. For both cases, the order of the scheme appears to be 1.
    
\begin{figure}[h!]
        \includegraphics[width=0.48\textwidth]{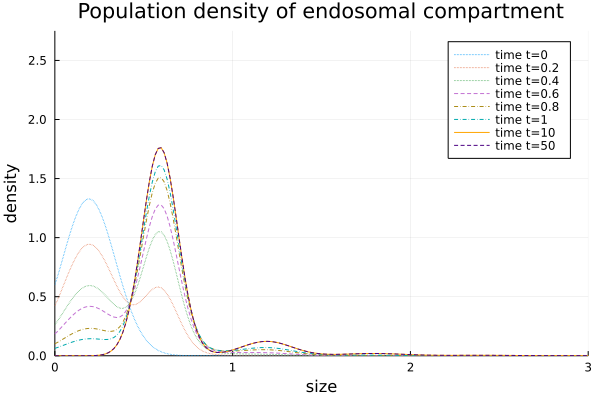}
	    \includegraphics[width=0.48\textwidth]{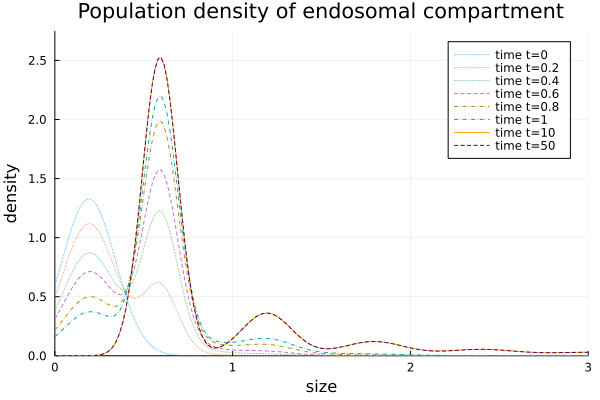}
    	\vspace{-0.5cm}
	    \caption{{\small\textit{ Evolution of the population density of endosomal compartement. On the left hand-side the condition Eq.~\eqref{hyp:1} is satisfied (we choose $\gamma \left(r\right) = \sqrt{3.1}$) whereas on the right hand-side the condition Eq.~\eqref{hyp:1} is not satisfied (we choose $\gamma\left(r\right) = 0.7$). The two graphics are obtained with the following parameters $\Delta t = 0.05$, $I^r=301$, $R=5$,  $\kappa\left(r,r'\right) = 1$, $\alpha\left(r,M\right) = M \cdot \mathcal{N}_r(0.6,0.1) $ and $f_0 \left(r\right) = 0.5 \cdot \mathcal{N}_r(0.2,0.15)$. In both cases, the curves for $t=10$ and $t=50$ are superimposed. }}}\label{fig:non_optimality_thm}
	    \vspace{-0.3cm}
        \end{figure}

    Finally, we test the optimality of the conditions from the long time behavior given in  Theorem \ref{theorem_longtime} for model \eqref{eq:coag_frag}. The detailed parametrization we use is given in the legend of Figure \ref{fig:non_optimality_thm}. First, we observe that we recover the convergence to a stationary state when the condition Eq.~\eqref{hyp:1} of Theorem \ref{theorem_longtime} is satisfied (left part of Figure \ref{fig:non_optimality_thm}). Numerically, it seems that the condition is too restrictive since the scheme stays stable in a wider range of parameters such that the condition Eq.~\eqref{hyp:1} is not satisfied. For example, on the right part of the Figure \ref{fig:non_optimality_thm} we have $3 \left\| \kappa \right\|_{L^\infty} \left\| \alpha  \right\|_{L^1} \geq \gamma_0^2 $, then  condition \eqref{hyp:1} is not satisfied but we still observe that the solution seems to converge to a stationary state.

\section{Applications}\label{sec:application}
    
     In this section we compare two different sets of experimental data with numerical simulations of our 2D model given by Eqs.~\eqref{eq:coag_frag_mol}-\eqref{eq:coag_frag_mol_membrane}. We obtain a good qualitative agreement. Our simulations also provide additional insight, which calls for new experiments.
    
    \subsection{Receptor trafficking}\label{traficking}

   The 2D model given by Eqs.~\eqref{eq:coag_frag_mol}-\eqref{eq:coag_frag_mol_membrane} is ideally suited to model receptor trafficking from plasma membrane to endocytic comparments. Briefly, the GPCRs are typically located at the surface of the cell, on the plasma membrane. Upon ligand binding, the GPCRs induce several signalling pathways as well as their own internalisation through clathrin-mediated endocytosis. The vesicles still carry the internalised receptors on their surface. Endosomes are then sorted thanks to complex processes which are not yet fully understood, but that depends both on the nature of the receptors and the ligand. Internalized receptors can indeed commit to several endosomal compartments of different kind, and be recycled at the cell surface, which could impact on the kinetic profile of the receptor. Consequently, endocytosis and post-endocytosis sorting regulate receptor cell surface density and signaling profile, and endosomal targeting of receptors may produce specificity in the signaling pathways. In particular, it has been shown \cite{jean-alphonseSpatiallyRestrictedProteincoupled2014} that the luteinizing hormone receptor (LHR) and the $\beta\text{2 -adrenergic}$ receptor (B2AR) are two GPCRs sorted to the regulated recycling pathway and undergo divergent trafficking to distinct endosomal compartments. B2AR traffics mostly to early endosomes (EEs) and LHR to pre-early endosomes (pre-EEs). The authors in \cite{jean-alphonseSpatiallyRestrictedProteincoupled2014} demonstrate that LHR endosome sizes increased over time quickly before reaching a plateau, producing a small endosome population (400-500 nm of diameter). The mean B2AR endosome sizes are bigger (1200-1400 nm of diameter). They also see that B2AR was more internalized than LHR in percentage but both receptors are equally recycled.\\

    In order to compare qualitatively our model with the experimental results presented in this paper, we choose the following parametrization (we recall the notation $\mathcal{N}_x(\mu,\sigma):x \mapsto\frac{1}{\sqrt{2 \pi }\sigma}e^{-\frac{1}{2}\left( \frac{x-\mu}{\sigma}\right)^2}$  and $\mathcal{P}_x(\bar{X},\epsilon):x \mapsto \left(\frac{\bar{X}+\epsilon-x}{\bar{X}+\epsilon}\right)^{\frac{1}{3}} \left(\frac{x}{\bar{X}}\right)^{\frac{2}{3}}$):
    \begin{itemize}
        \item Endosomes are created with a size following a Gaussian law with mean $\mu=200 nm^3$ and standard deviation $\sigma=10 nm^3$, and with a quantity of reactant proportional to their size. That is, $\alpha(r,a,M)=\overline{\alpha}\mathcal{N}_r(200,10)\times\mathcal{N}_a(r,0.5)\times M$ for some rate $\overline{\alpha}>0$ that will depend on the receptor.
        \item Endosome recycling increases with the surface of endosomes and decreases with their volume ( $\lambda \propto r^{\frac{2}{3}}-r$), e.g. $\lambda :=10^{-2} \times \mathcal{P}_r(2000,0) $.
        \item Endosomes fuse at a constant coagulation kernel $\kappa$, whose rate will depend on the receptor.
        \item Small endosomes are degraded at constant (low) rate, and degradation rates quickly increases for larger endosomes, e.g. $\gamma:=10^{-5}+2\times 10^{1}\times \left(\frac{r-1950}{50}\right)^4 \indic{\left\{r>1950\right\}}$.
        \item We don't consider reactions, e.g. $V$ and $J_M$ are taken as null functions. 
        \end{itemize}
 
    We choose as initial conditions $f_0=0$ and $M_0=7.2\times 10^{-4}$. For the numerical scheme, we take $\Delta t= 10^{-1}$, $I^r=30$, $I^a=30$, $R= 2000$, $A=2000$.
        
\begin{figure}[h!]
	\includegraphics[width=0.45\linewidth]{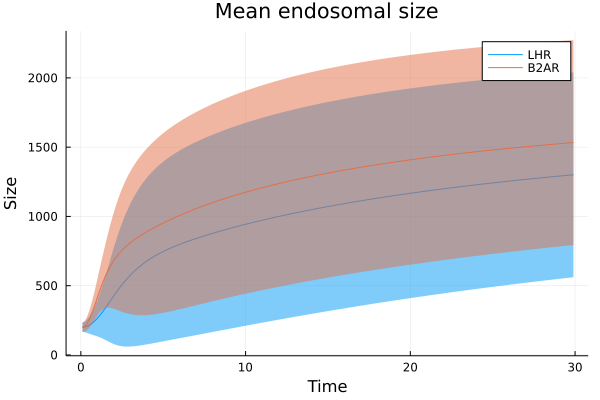}
	\includegraphics[width=0.45\linewidth]{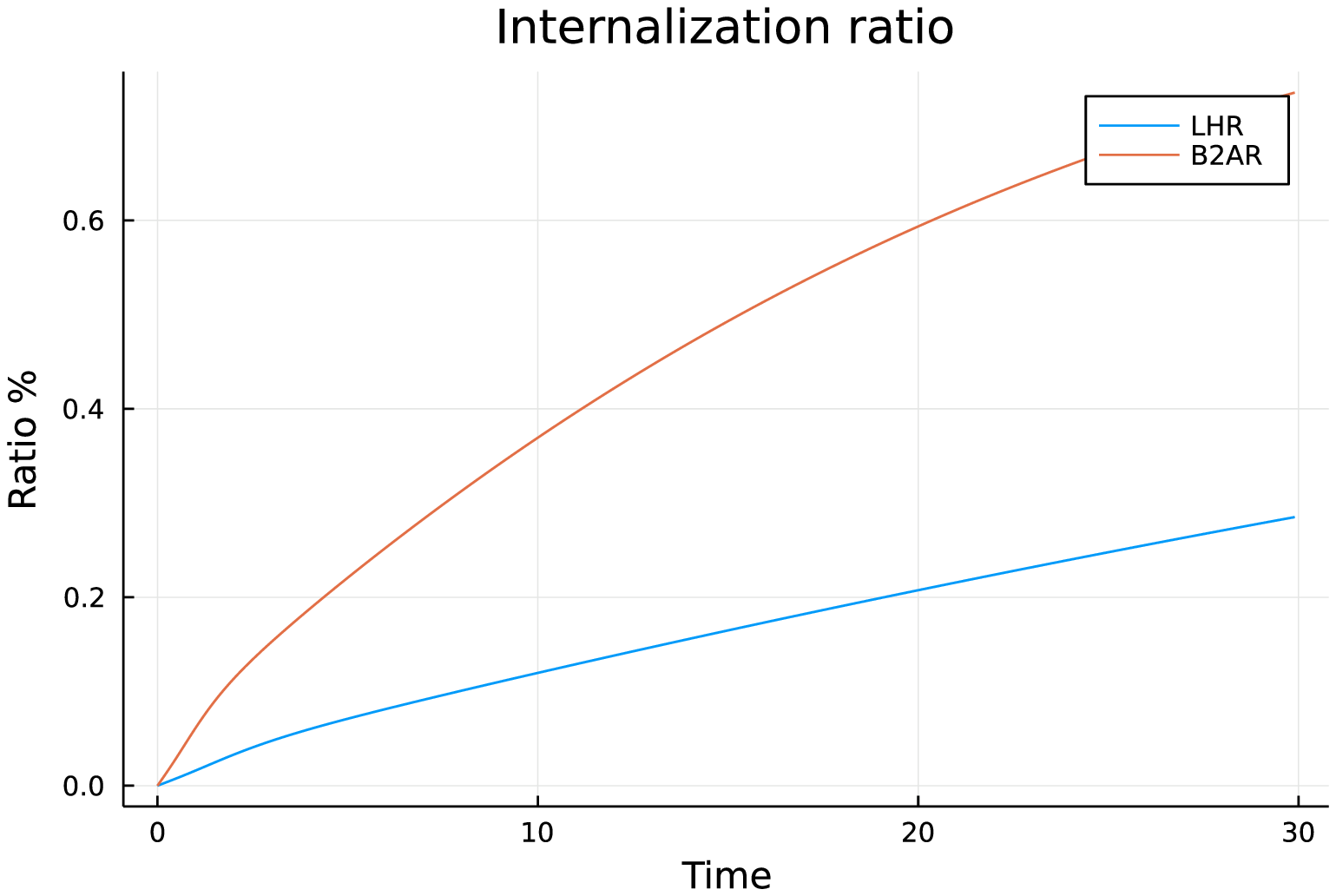}
	\vspace{-0.5cm}
	\caption{{\small\textit{ LHR and B2AR affect only the internalisation rate - The mean endosomal size corresponds to  $\left\|f\right\|_{1r}(t)$, normalised by the total mass $\left\|f\right\|_{0}(t)$ (the standard deviation of the size is represented in light colours). The internalisation ratio corresponds to $\frac{\left\|f\right\|_{1a}(t)}{M_0}$ (we suppose that the process of exocytosis is negligible at this timescale).}}}\label{fig1}
	\vspace{-0.3cm}
    \end{figure}

\begin{figure}[h!]
	\includegraphics[width=0.45\textwidth]{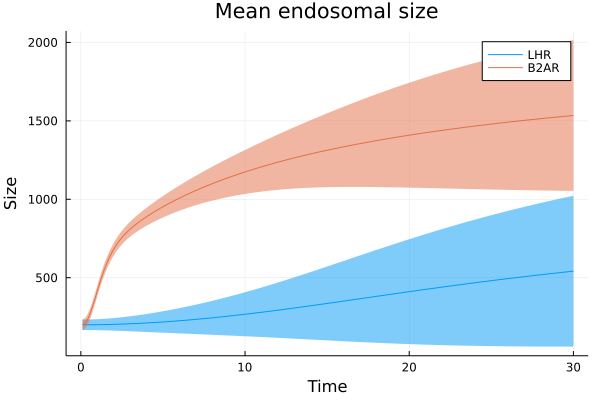}
	\includegraphics[width=0.45\textwidth]{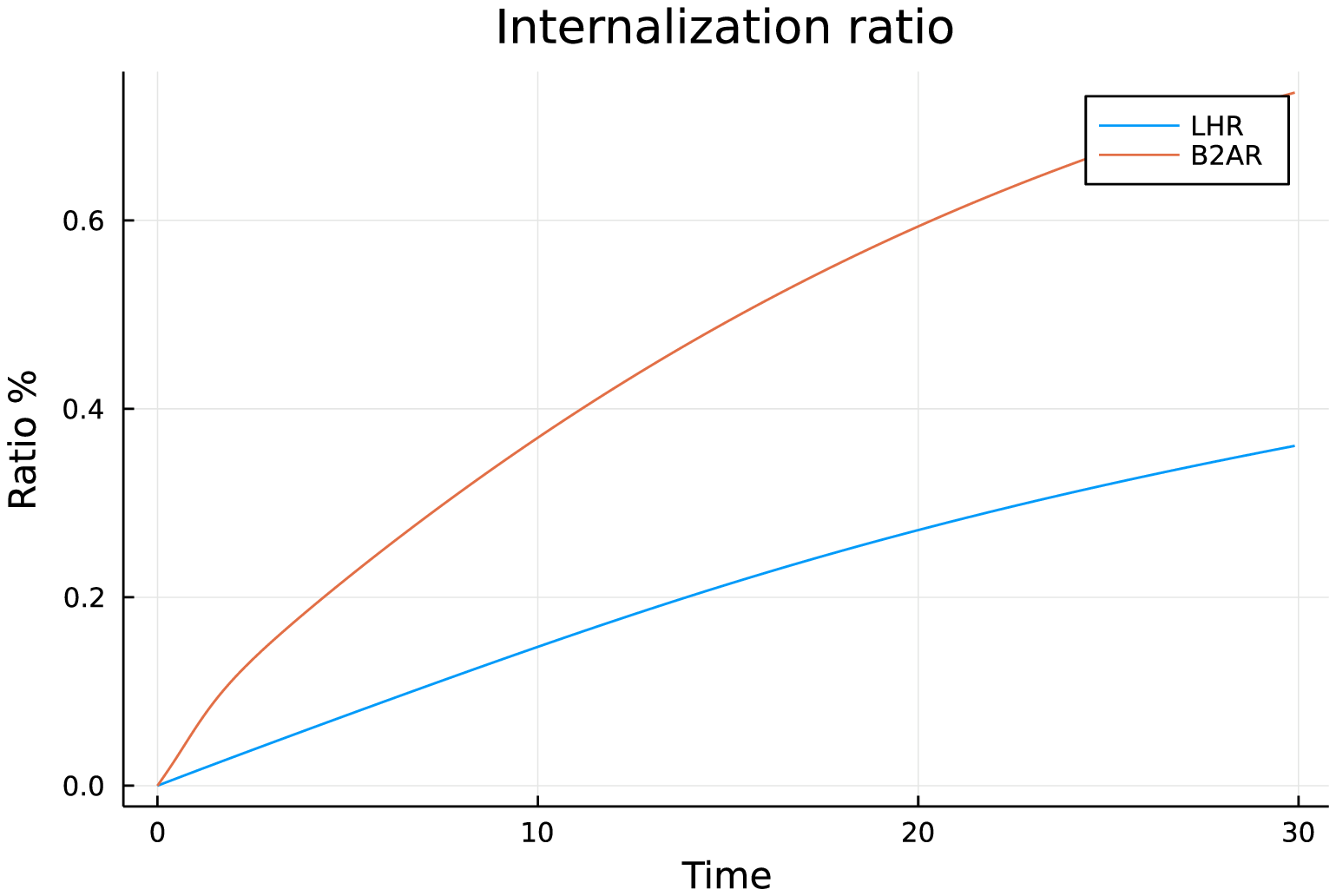}
	\vspace{-0.5cm}
	\caption{{\small\textit{ LHR and B2AR affect the internalisation and the coagulation rate - The mean endosomal size corresponds to   $\left\|f\right\|_{1r}(t)$, normalised by the total mass $\left\|f\right\|_{0}(t)$ (the standard deviation of the size is represented in light colours). The internalisation ratio correspond to $\frac{\left\|f\right\|_{1a}(t)}{M_0}$ (we suppose that the process of exocytosis is negligible at this timescale).}}}\label{fig2}
	\vspace{-0.3cm}
    \end{figure}

    To explain the qualitative differences between LHR and B2AR trafficking observed in \cite{jean-alphonseSpatiallyRestrictedProteincoupled2014}, we tested two hypotheses here. 
    \begin{itemize}
        \item (H1) In the first hypothesis (Figure \ref{fig1}), we modify only the internalisation rate $\overline{\alpha}>0$ between LHR and B2AR, with a higher internalisation rate for the B2AR, keeping all the remaining parameters the same. This hypothesis is in line with \cite{jean-alphonseSpatiallyRestrictedProteincoupled2014}. LHR and B2AR affect only the endocytosis rate as follows:
        \newline
        {\centering
        \begin{tabular}{|c|c|c|}
            \hline
            Parameters & $\kappa$ & $\overline{\alpha}$  \\
            \hline
            LHR & $5\times 10^{-1}$ & $8\times 10^{-5} $ \\
            \hline
            B2AR & $5\times 10^{-1}$ & $3\times 10^{-4} $ \\
            \hline
        \end{tabular}
        }
        \item (H2) In the second hypothesis (Figure \ref{fig2}), both the internalisation rate and the coagulation rate are higher for the B2AR compared to the LHR (H2). LHR and B2AR affect the internalisation and the coagulation rate as follows:  
        \newline
        {\centering
        \begin{tabular}{|c|c|c|}
            \hline
            Parameters & $\kappa$ & $\overline{\alpha}$\\
            \hline
            LHR & $5\times 10^{-3}$ & $8\times 10^{-5}$\\
            \hline
            B2AR & $5\times 10^{-1}$ & $3\times 10^{-4}$\\
            \hline
        \end{tabular}
        }
    \end{itemize}
    In both hypotheses the internalization ratio presented in \cite{jean-alphonseSpatiallyRestrictedProteincoupled2014} is well reproduced by the model. The major discrepancy of the first hypothesis with the experiments is the variance in size of the endosomal population, which seems too high compared to experimentation. Furthermore the production of large endosomes with B2AR is not high enough (Figure \ref{fig1}). In the second hypothesis, however, these two  discrepancies are not present anymore (Figure \ref{fig2}). These results indicate that the differences in endosome dynamics between LHR and B2AR signaling pathways seem not only due to a difference in internalisation rate, but probably also to the coagulation dynamics inside cells. Difference in coagulation dynamics may be explained by differences in endosomes sorting and/or differences in molecular composition of endosomes, which are believed to be of different nature between LHR and B2AR vesicles \cite{jean-alphonseSpatiallyRestrictedProteincoupled2014}.\\

        \subsection{Second effector signaling}\label{signaling}
        %D.White paper and Vilardaga paper
        The second application of our model we present here concerns
        %is
        the efficacy of second messenger molecules as a function of the localization of the active receptor. We focus on the production of cyclic adenosine monophosphate (cAMP) induced by the activation of the parathyroid hormone receptor (PTHr). Recent discoveries find that the PTHr may engage cAMP signaling not only at the cell’s plasma membrane but also in early endosomes after receptor internalisation through endocytosis \cite{vilardagaEndosomalGenerationCAMP2014}.
Furthermore, the full-length PTH (${PTH}^{WT}$ or LA-PTH) induces through the activation of PTHr an augmentation of production of cAMP in the endosomes, whereas the ligand ${PTH}^{7d}$ (PTHR peptide ligand through amino acid epimerization at position 7 of ${PTH}^ {1-34}$) induces the production of cAMP at the plasma membrane.  In White’s paper \cite{whiteSpatialBiasCAMP2021}, the authors show that even if the production place is different, the total amount of cAMP stays the same after a small time, a phenomenon that could be named location-biased, and that can have implications for the cellular response.   

We  could reproduce qualitatively these observations with our model following  this parametrization:
    \begin{itemize}
        \item Endosomes are created with a size following a Gaussian law with mean $\mu = 200 nm^3$ and standard deviation $\sigma=10 nm^3$, and with a quantity of reactant following the positive part of a Gaussian law around 0 with standard deviation $\sigma=0.1$ (arbitrary unit), e.g. $\alpha:=\mathcal{N}_r(200,10)\times\mathcal{N}_a(0,0.1))\times M$
        \item Endosome recycling and degradation increase with the surface of endosomes and decrease with their volume ($\lambda \text{ and } \gamma \propto r^{\frac{2}{3}}-r$), e.g. $\lambda :=10^{-2}\times \mathcal{P}_r(2000,100) $,  and $\gamma:= \mathcal{P}_r(2000,100) \indic{\left\{r\leq 1950\right\}}+2\times 10^{2}\times \left(\frac{r-1950}{50}\right)^4 \indic{\left\{r>1950\right\}}$.
        
        \item cAMP is produced at the plasma membrane at constant rate $v_M$ and linearly degraded, with a saturation of $\bar{M}$, e.g $J_M(M)= v_M\times \left(\frac{\bar{M}-M}{\bar{M}}\right)$. Both $\bar{M}$, and $v_M$ may depend on the ligand.
        \item cAMP is produced in endosomes at two different rates, $v_s$, for the endosomes smaller than $\bar{r}$ and $v_l$, for larger endosomes. Also the amount of saturation depends linearly on the size of the endosomes (caracterised by a portion $p$ of r). Hence $V(r,a):=\left( v_s\indic{\left\{\epsilon\leq r\leq\bar{r}\right\}} + v_l\indic{\left\{\bar{r}<r\right\}}\right)\times \left(1-\frac{a}{pr}\right)$ with $0<\epsilon\ll1$. Here $v_s$, $v_l$ and $p$ may depend on the ligand.
        \item Endosomes fuse via a constant coagulation kernel $\kappa:= 2\times 10^{-1}$. 

    \end{itemize}
    
We start from a zero initial condition, $f_0$ is the null function and $M_0=0$. For the numerical scheme, we take $\Delta t=3 \times 10^{-2}$, $I^r=30$, $I^a=30$, $R= 2000$, $A=30$.

\begin{figure}[h!]
    \includegraphics[width=0.30\textwidth]{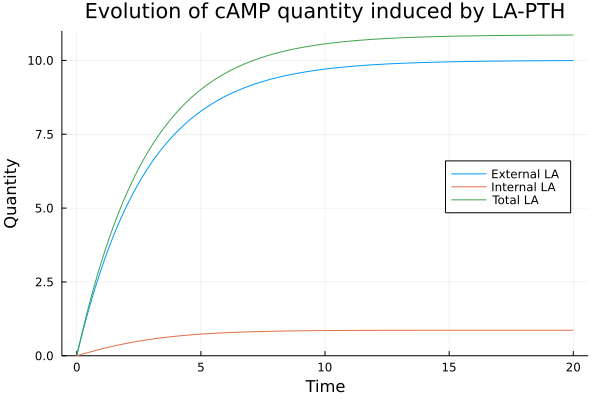}
    \includegraphics[width=0.30\textwidth]{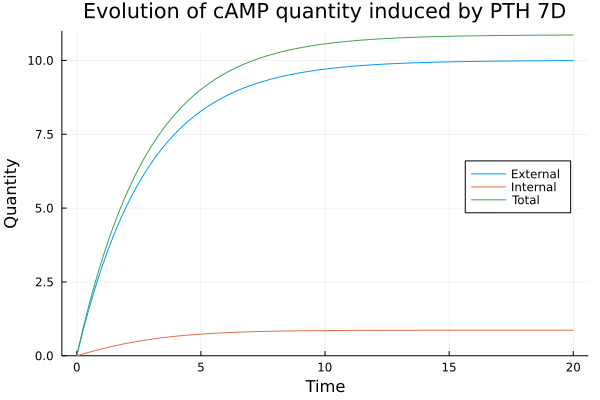}
    \includegraphics[width=0.30\textwidth]{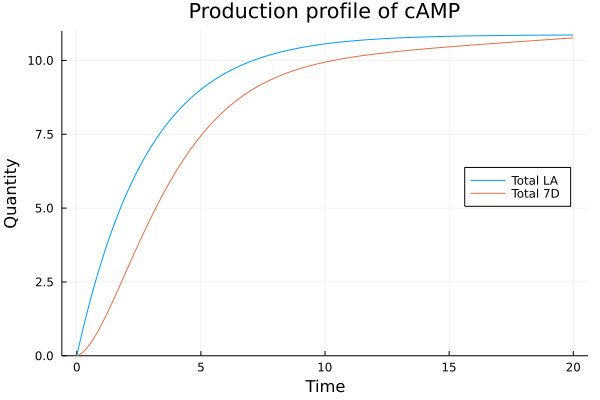}
	\vspace{-0.5cm}
    \caption{{\small\textit{ Comparison in the production of cAMP depending on the production area - PTH 7D and LA-PTH affect only the production rate - The internal quantity of cAMP is given by $\left\|f\right\|_{1a}(t)$ and the external one by $M$.}}}\label{same_satur}
    \vspace{-0.3cm}
\end{figure}

\begin{figure}[h!]
    \includegraphics[width=0.30\textwidth]{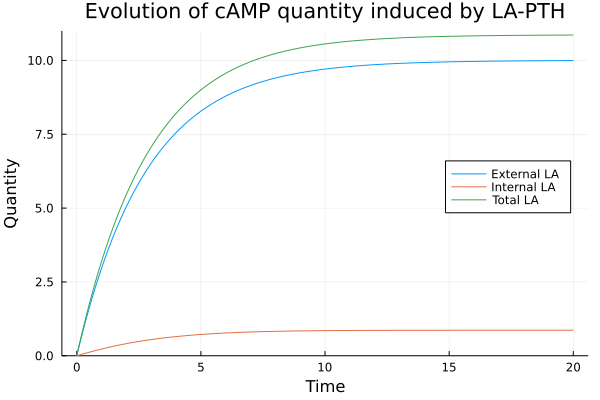}
    \includegraphics[width=0.30\textwidth]{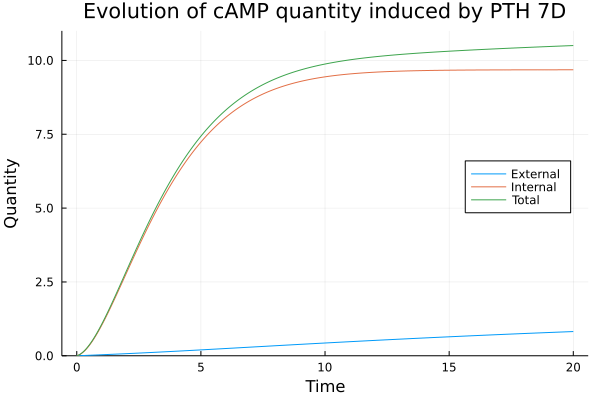}
    \includegraphics[width=0.30\textwidth]{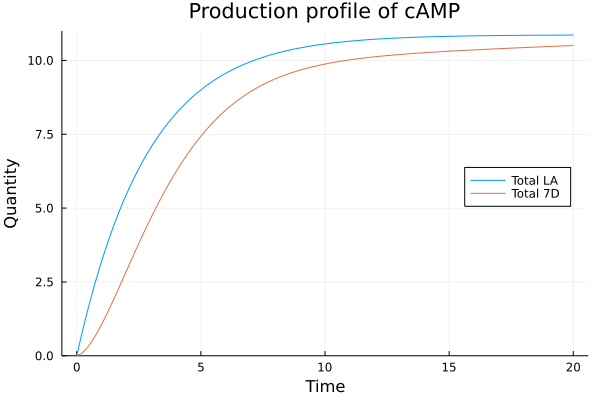}
	\vspace{-0.5cm}
    \caption{{\small\textit{ Comparison in the production of cAMP depending on the production area - PTH 7D and LA-PTH affect both the production in term of rate and capacity of production  - The internal quantity of cAMP is given by $\left\|f\right\|_{1a}(t)$ and the external one by $M$.}}}\label{diff_satur}
	\vspace{-0.3cm}
\end{figure}

We have two different hypotheses to explain the behaviour described above:
    \begin{itemize}
        \item (H1) Suppose LA-PTH and PTH 7D differs in cAMP production kinetics only in terms of rate (Figure \ref{same_satur}), with a higher rate for LA-PTH at the plasma membrane ($v_M^{PTH 7D}<v_M^{LA-PTH}$), and a higher rate for PTH 7D in the endosomes ($v_s^{LA-PTH}<v_s^{PTH 7D}$ and $v_l^{LA-PTH}<v_l^{PTH 7D}$).
        \newline
                {\centering
        \begin{tabular}{|c|c|c|c|c|c|}
            \hline
            Parameters &   $v_s$ & $v_l$ & $p$ & $v_M$ &  $\bar{M}$\\
            \hline
            LA-PTH &   $0.05$ & $0.02$ & $1/20$ & $3.5$ &  $10$\\
            \hline
            PTH 7D &   $5$ & $2$ & $1/20$ & $0.035$ &  $10$\\
            \hline
    	    \end{tabular}
        }
        
        \item (H2) Suppose LA-PTH and PTH 7D differs in cAMP production kinetics not only in terms of rate but also in terms of saturation (Figure \ref{diff_satur}) with $\bar{M}^{LA-PTH}>\bar{M}^{{PTH}^{7D}}$ and $p^{{PTH}^{7D}}>p^{LA-PTH}$ .
                \newline
                {\centering
        \begin{tabular}{|c|c|c|c|c|c|}
            \hline
            Parameters &   $v_s$ & $v_l$ & $p$ & $v_M$ &  $\bar{M}$\\
            \hline
            LA-PTH &   $0.5$ & $0.2$ & $1/200$ & $3.5$ &  $10$\\
            \hline
            PTH 7D &   $5$ & $2$ & $1/20$ & $0.35$ &  $1$\\
            \hline
    	    \end{tabular}
        }
    \end{itemize}

With both hypotheses, we observe a much more efficient cAMP production at the plasma membrane with LA-PTH and a much more efficient cAMP production in the endosomes with PTH 7D (Figures \ref{same_satur} and \ref{diff_satur}). Consistently with the observation in \cite{whiteSpatialBiasCAMP2021}, both total responses have similar magnitude for the time period of the numerical simulation. 

However, from the numerical simulation presented in figures \ref{same_satur} and \ref{diff_satur}, the cAMP production has already reached a "stable" state at $T=20$ for LA-PTH, while it keeps increaseasing for PTH 7D. Therefore, a longer time measurement could discriminate between both ligands.

Also we could notice a fine kinetic difference between the responses induces by LA-PTH and PTH 7D with the two hypotheses. Indeed, PTH 7D leads to a convex kinetic production of cAMP during the early dynamics, which switches to a concave kinetic at later time. Whereas with LA-PTH the production stays concave all time long. Of course this behaviour may be quite complicated to observe experimentally due to the  accuracy of the measures.

\section{Discussion}

%%-----------------------------
%%      your bibliography
%%-----------------------------
\bibliographystyle{ieeetr}	
%\bibliography{projet_AE_bold.bib}
\bibliography{projet_AE.bib}

\begin{thebibliography}{10}

\bibitem{kenakinBiasedReceptorSignaling2019}
T.~Kenakin, ``Biased {Receptor} {Signaling} in {Drug} {Discovery},'' {\em
  Pharmacological Reviews}, vol.~71, pp.~267--315, Apr. 2019.

\bibitem{jean-alphonseRegulationGPCRSignal2011}
F.~{Jean-Alphonse} and A.~C. Hanyaloglu, ``Regulation of {{GPCR}} signal
  networks via membrane trafficking,'' {\em Molecular and Cellular
  Endocrinology}, vol.~331, pp.~205--214, Jan. 2011.

\bibitem{vilardagaEndosomalGenerationCAMP2014}
J.-P. Vilardaga, F.~{Jean-Alphonse}, and T.~J. Gardella, ``Endosomal generation
  of {{cAMP}} in {{GPCR}} signaling,'' {\em Nat Chem Biol}, vol.~10,
  pp.~700--706, Sept. 2014.

\bibitem{kimSpatiotemporalCharacterizationGPCR2021}
H.~Kim, H.~N. Lee, J.~Choi, and J.~Seong, ``Spatiotemporal {{Characterization}}
  of {{GPCR Activity}} and {{Function}} during {{Endosomal Trafficking
  Pathway}},'' {\em Anal. Chem.}, p.~acs.analchem.0c03323, Jan. 2021.

\bibitem{jean-alphonseSpatiallyRestrictedProteincoupled2014}
F.~{Jean-Alphonse}, S.~Bowersox, S.~Chen, G.~Beard, M.~A. Puthenveedu, and
  A.~C. Hanyaloglu, ``Spatially {{Restricted G Protein-coupled Receptor
  Activity}} via {{Divergent Endocytic Compartments}},'' {\em J Biol Chem},
  vol.~289, pp.~3960--3977, Feb. 2014.

\bibitem{lygaPersistentCAMPSignaling2016}
S.~Lyga, S.~Volpe, R.~C. Werthmann, K.~G{\"o}tz, T.~Sungkaworn, M.~J. Lohse,
  and D.~Calebiro, ``Persistent {{cAMP Signaling}} by {{Internalized LH
  Receptors}} in {{Ovarian Follicles}},'' {\em Endocrinology}, vol.~157,
  pp.~1613--1621, Apr. 2016.

\bibitem{sayersIntracellularFollicleStimulatingHormone2018}
N.~Sayers and A.~C. Hanyaloglu, ``Intracellular {{Follicle-Stimulating Hormone
  Receptor Trafficking}} and {{Signaling}},'' {\em Frontiers in Endocrinology},
  vol.~9, Nov. 2018.

\bibitem{whiteSpatialBiasCAMP2021}
A.~D. White, K.~A. Pe{\~n}a, L.~J. Clark, C.~S. Maria, S.~Liu, F.~G.
  {Jean-Alphonse}, J.~Y. Lee, S.~Lei, Z.~Cheng, C.-L. Tu, F.~Fang, N.~Szeto,
  T.~J. Gardella, K.~Xiao, S.~H. Gellman, I.~Bahar, I.~Sutkeviciute, W.~Chang,
  and J.-P. Vilardaga, ``Spatial bias in {{cAMP}} generation determines
  biological responses to {{PTH}} type 1 receptor activation,'' {\em Science
  Signaling}, vol.~14, p.~eabc5944, Oct. 2021.

\bibitem{sorkinEndocytosisSignallingIntertwining2009}
A.~Sorkin and M.~{von Zastrow}, ``Endocytosis and signalling: Intertwining
  molecular networks,'' {\em Nat Rev Mol Cell Biol}, vol.~10, pp.~609--622,
  Sept. 2009.

\bibitem{birtwistleEndocytosisSignallingMeeting2009}
M.~R. Birtwistle and B.~N. Kholodenko, ``Endocytosis and signalling: {{A}}
  meeting with mathematics,'' {\em Mol Oncol}, vol.~3, pp.~308--320, Aug. 2009.

\bibitem{villasenorSignalProcessingEndosomal2016}
R.~Villase{\~n}or, Y.~Kalaidzidis, and M.~Zerial, ``Signal processing by the
  endosomal system,'' {\em Current Opinion in Cell Biology}, vol.~39,
  pp.~53--60, Apr. 2016.

\bibitem{ingallsMathematicalModelingSystems2013}
B.~Ingalls, {\em Mathematical {Modeling} in {Systems} {Biology}}.
\newblock MIT Press, 2013.

\bibitem{feinbergFoundationsChemicalReaction2019}
M.~Feinberg, {\em Foundations of {Chemical} {Reaction} {Network} {Theory}}.
\newblock Applied {Mathematical} {Sciences}, Springer International Publishing,
  2019.

\bibitem{andersonStochasticAnalysisBiochemical2015}
D.~F. Anderson and T.~G. Kurtz, {\em Stochastic {Analysis} of {Biochemical}
  {Systems}}.
\newblock No.~1.2 in Mathematical {Biosciences} {Institute} {Lecture} {Series},
  Springer International Publishing, 2015.

\bibitem{kholodenkoCellsignallingDynamicsTime2006}
B.~N. Kholodenko, ``Cell-signalling dynamics in time and space,'' {\em Nature
  Reviews Molecular Cell Biology}, vol.~7, pp.~165--176, Mar. 2006.

\bibitem{weddellIntegrativeMetamodelingIdentifies2017}
J.~C. Weddell and P.~I. Imoukhuede, ``Integrative meta-modeling identifies
  endocytic vesicles, late endosome and the nucleus as the cellular
  compartments primarily directing {{RTK}} signaling,'' {\em Integrative
  Biology}, vol.~9, pp.~464--484, May 2017.

\bibitem{foretGeneralTheoreticalFramework2012}
L.~Foret, J.~E. Dawson, R.~Villase{\~n}or, C.~Collinet, A.~Deutsch, L.~Brusch,
  M.~Zerial, Y.~Kalaidzidis, and F.~J{\"u}licher, ``A general theoretical
  framework to infer endosomal network dynamics from quantitative image
  analysis,'' {\em Curr Biol}, vol.~22, pp.~1381--1390, Aug. 2012.

\bibitem{dusoStochasticReactionNetworks2020}
L.~Duso and C.~Zechner, ``Stochastic reaction networks in dynamic compartment
  populations,'' {\em Proceedings of the National Academy of Sciences},
  vol.~117, pp.~22674--22683, Sept. 2020.

\bibitem{pietzschCompartorToolboxAutomatic2021}
T.~Pietzsch, L.~Duso, and C.~Zechner, ``Compartor: a toolbox for the automatic
  generation of moment equations for dynamic compartment populations,'' {\em
  Bioinformatics}, vol.~37, pp.~2782--2784, Sept. 2021.

\bibitem{alexandrovDynamicsIntracellularClusters2022}
D.~V. Alexandrov, N.~Korabel, F.~Currell, and S.~Fedotov, ``Dynamics of
  intracellular clusters of nanoparticles,'' {\em Cancer Nanotechnology},
  vol.~13, p.~15, May 2022.

\bibitem{chae_existence_1995}
D.~Chae and P.~B. Dubovskii, ``Existence and uniqueness for spatially
  inhomogeneous coagulation equation with sources and effluxes,'' {\em
  Zeitschrift f\"ur angewandte Mathematik und Physik {ZAMP}}, vol.~46, no.~4,
  pp.~580--594, 1995-07-01.

\bibitem{gajewski_first_1983}
H.~Gajewski, ``On a first order partial differential equation with nonlocal
  nonlinearity,'' {\em Mathematische Nachrichten}, vol.~111, no.~1,
  pp.~289--300, 1983.

\bibitem{dubovskii_mathematical_nodate}
P.~B. Dubovskii, ``Mathematical theory of coagulation,'' {\em Lectures notes
  from National Univ}, 1993-1994.

\bibitem{laurencot_stationary_2020}
P.~Laurençot, ``Stationary solutions to smoluchowski's coagulation equation
  with source,'' {\em North-Western European Journal of Mathematics}, vol.~6,
  p.~137, 2020.

\bibitem{Ofelia2015Thesis}
M.~O. V\'asquez, {\em Ecuaciones en derivadas parciales para el an\'alisis de
  modelos biopolim\'ericos}.
\newblock PhD thesis, Universidad de Granada, 2015.

\bibitem{ghosh_equilibrium_2023}
D.~Ghosh, J.~Paul, and J.~Kumar, ``On equilibrium solution to a singular
  coagulation equation with source and efflux,'' {\em Journal of Computational
  and Applied Mathematics}, vol.~422, p.~114909, 2023-04-01.

\bibitem{collet_lifshitz-slyozov_1999}
J.-F. Collet and T.~Goudon, ``Lifshitz-slyozov equations: The model with
  encounters,'' {\em Transport Theory and Statistical Physics}, vol.~28, no.~6,
  pp.~545--573, 1999-10-01.

\bibitem{Canizo2006Thesis}
J.~A. Canizo, {\em Some problems related to the study of interaction kernels:
  coagulation, fragmentation and diffusion in kinetic and quantum equations}.
\newblock PhD thesis, Universidad de Granada, 2006.

\bibitem{bourgadeConvergenceFiniteVolume2007}
J.-P. Bourgade and F.~Filbet, ``Convergence of a finite volume scheme for
  coagulation-fragmentation equations,'' {\em Mathematics of Computation},
  vol.~77, pp.~851--883, Dec. 2007.

\bibitem{hingantDerivationMathematicalStudy2015}
E.~Hingant and M.~Sepúlveda, ``Derivation and mathematical study of a
  sorption-coagulation equation,'' {\em Nonlinearity}, vol.~28, pp.~3623--3661,
  Oct. 2015.

\end{thebibliography}

\end{document}